\documentclass{article}

\usepackage[utf8]{inputenc} 
\usepackage{amssymb,amsmath} 
\usepackage{graphicx} 
\usepackage{array} 
\usepackage{eurosym} 
\usepackage{xcolor} 
\usepackage{booktabs}
\usepackage{epstopdf}
\usepackage{pifont} 
\usepackage{algorithm} 
\usepackage{algpseudocode}
\usepackage{multicol}
\usepackage{multirow}
\usepackage{subcaption}
\usepackage{changepage}
\usepackage{hyperref}
\usepackage{adjustbox}
\usepackage{graphicx}
\usepackage{multirow}
\usepackage{amsmath} 
\usepackage{textcomp}
\usepackage[affil-it]{authblk}

\begin{document}

\title{Operational planning and bidding for district heating systems with uncertain renewable energy production}

\author[1]{Ignacio~Blanco}
\author[1,*]{Daniela Guericke}
\author[2]{Anders N. Andersen}
\author[1]{Henrik~Madsen}

\affil[1]{\small Technical University of Denmark, Department for Applied Mathematics and Computer Science, Richard Petersens Plads, 2800 Kgs. Lyngby, Denmark}
\affil[2]{\small EMD International A/S, Niels Jernesvej 10, 9220 Aalborg \O, Denmark}
\affil[*]{\small Corresponding author: Daniela Guericke, dngk@dtu.dk}

\maketitle

\begin{abstract}
In countries with an extended use of district heating (DH), the integrated operation of DH and power systems can increase the flexibility of the power system achieving a higher integration of renewable energy sources (RES). DH operators can not only provide flexibility to the power system by acting on the electricity market, but also profit from the situation to lower the overall system cost. However, the operational planning and bidding includes several uncertain components at the time of planning: electricity prices as well as heat and power production from RES. In this publication, we propose a planning method that supports DH operators by scheduling the production and creating bids for the day-ahead and balancing electricity markets. The method is based on stochastic programming and extends bidding strategies for virtual power plants to the DH application. The uncertain factors are considered explicitly through scenario generation. We apply our solution approach to a real case study in Denmark and perform an extensive analysis of the production and trading behaviour of the DH system. The analysis provides insights on how DH system can provide regulating power as well as the impact of uncertainties and renewable sources on the planning. Furthermore, the case study shows the benefit in terms of cost reductions from considering a portfolio of units and both markets to adapt to RES production and market states.

\end{abstract}

\section{Introduction}

To achieve the decarbonization of the energy sector, several countries especially in the European Union started to consider district heating (DH) and cooling systems for CO2-emissions reduction strategies \cite{studyone27:online}. Since it is assumed that fossil fuels will be mostly replaced by intermittent renewable energy sources (RES), DH and cooling systems can facilitate a larger share of intermittent energy sources in the energy mix following the concept of integrated energy systems \cite{CONNOLLY2014475}. DH systems are able to contribute to the grid balancing by the use of flexible heat and power production, power-to-heat technologies and thermal storages. 

The efficiency of DH systems has been demonstrated already in countries in northern Europe.  In Denmark, more than 60\% of the heat consumption is delivered by DH \cite{euroheat} and there exist a total of approximately 400 DH systems. The major part of those are small/medium DH systems that are usually operated based on a portfolio of different units such as CHP units (e.g. gas engines), fuel boilers, and power-to-heat technologies such as electric boilers and heat pumps. Also the installation of large solar thermal facilities ($\geq$1000$\text{ m}^{2}$) in Denmark has increased significantly during the last years and it is expected that 20\% of the total heat consumption will be covered by solar heating in 2025 \cite{regulati24:online}. Furthermore, the wider spread of power-to-heat technologies and decentralization of power production enables DH providers to include renewable power production, e.g., in the form of wind farms, to their portfolio. Although the primary goal of the DH operator is to fulfill the heat demand in the DH network at lowest cost, selling the power production from the CHP units or other RES as well as buying the power for heat-to-power technologies on electricity markets offers the potential for additional income resulting in lower total operating costs. 
However, as the RES production in the power and heat systems depends on weather conditions, the operation and planning has to deal with an increased complexity and uncertainty, which requires advanced modeling techniques \cite{madsen2015integrated}.

In this publication, we pursue two main objectives. First, we propose an operational planning method for DH operators coping with the complexity of a system with several traditional and RES production units. This includes the bidding in two electricity markets, namely day-ahead and balancing market. The method uses stochastic programming to capture the uncertainties and is based on models proposed for virtual power plants (VPPs) \cite{pandvzic2013offering}. Second, we use the proposed method to analyze the real case of a district heating system in  Hvide Sande, Denmark. The analysis investigates among others the behaviour of the DH system in different situations, the influence of uncertainty in the RES production and benefits from including RES power production. The results offer several insights on how DH systems should operate and can benefit in future systems with high shares of RES .

\subsection{Description of electricity markets}

Nowadays, the integration of the power and DH system is achieved through the participation of the latter in the electricity markets. Before describing the related work, we want to recall the concepts of the day-ahead and balancing electricity markets that are considered by the proposed planning method. 

In most of the EU countries, the short-term trading of electricity is organized in a similar way. Most of the power volume is traded one day before the energy is delivered in the so-called \textit{day-ahead market}. To ensure enough backup generation, producers can also bid offers in the \textit{reserve capacity market} which takes place usually also one day before the delivery of energy. Getting closer to the time of delivery, \textit{intra-day markets} are organized throughout the day to help RES producers submit more accurate power production offers. The purpose of these markets is to correct the imbalances produced by RES allowing producers to reformulate their bids. Finally, \textit{balancing markets} are organized each hour of the day with gate closure one hour before the energy delivery.

\textit{Balancing markets} are slightly different from intra-day markets and take place shortly before hour of energy delivery. The balancing markets are cleared by the TSO and their goal is to provide flexibility for the operation of the system and not to the producers as it is the case for the intra-day market. Balancing prices are highly volatile and quite unpredictable even for the following hour. In addition, balancing offers are just activated in case the TSO has need for regulation. In the case that there is a lack of power production due to a failure of a unit or an unpredictable  demand, the TSO will activate offers for upward regulation paying producers to increase the production of their power plants. On the contrary, if there is more power production than expected due to an excess of RES production, the TSO will activate downward regulation offers for producers to deactivate the production they had previously scheduled in the day-ahead market or incentive more power consumption. 

To efficiently operate in these markets, producers and consumers are allowed to submit price dependent bids. These type of bids consist of pair-wise points of power volume and power prices that must follow a merit ascending or descending order. In this way, producers and consumers are able to provide a wider range of offers to hedge against the uncertain electricity prices. In this work, we focus on day-ahead and balancing markets.

\subsection{Related work}

As mentioned before, the integration of heat and power production units complicates the operation of the system requiring suitable tools. Among other techniques, mixed integer linear programming (MILP) has been shown as one well-suited approach to optimize the operation of DH systems. To provide some examples, the authors in \cite{carpaneto2015optimal} propose a unit commitment model that optimizes the integration of a solar collector in a DH system that includes one fuel boiler and one CHP unit connected to a thermal storage tank. Furthermore, the authors in \cite{wang2015modelling} work with a DH system that includes several CHP units, fuel boilers, thermal storage as well as a solar thermal plant. The authors propose an optimization model that accounts for the synchronization of the operation of the units providing an extensive analysis of flexibility between units. Finally, the authors in \cite{li2016optimal} go a step further by introducing a wind farm in a DH system that can feed both a heat pump and the power grid. To provide flexibility to the system, they also integrate a CHP unit and a thermal storage that increases the complexity of operating the system. All these presented publications have in common that they operate a portfolio of distributed generators and flexible loads.

Apart from the operational planning, planning methods have to consider the bidding in electricity markets. Nowadays, producers often base their offers on the given electricity price forecast, which is very volatile due to the variability of RES production and uncertain one day before the energy is delivered \cite{paraschiv2014impact}. Additionally, the production from RES in the DH system itself is uncertain. Consequently, tools that optimize the operation of DH systems and propose bidding strategies need to consider the uncertainty given by price and production. Despite several bidding strategies for price-taker power producers in the day-ahead market have been proposed, (see references in \cite{Kwon2012}), the authors in \cite{pandvzic2013offering} demonstrate that under high uncertainty of electricity prices the use of stochastic programming \cite{birge2011introduction} for creating bidding curves for the day-ahead market renders good solutions that consider the uncertainty involved in the bidding process. Based on the representation of the uncertain electricity prices as scenarios, the authors use non-anticipativity constraints that order the bids presented to the market in a step-wise manner to create price dependent bids.

The above mentioned methods consider power production only. Hence, they are not directly applicable for DH operators as the heat production is neglected. The heat production is an important part and a planning method needs to ensure heat demand fulfillment as well as consider the limitations of the production units and storages. Therefore, bidding methods for systems with a connected DH system need to model the heat production as well. For example, the method proposed in \cite{schulz2016optimal} determines the optimal production of a CHP unit. The bidding price is the price forecast, which is the same price used to determined the power production. In \cite{ravn2004modelling} the authors propose a bidding strategy for CHP units that takes into account other heat units to define the heat production costs to determine the bidding price. Finally, the authors in \cite{dimoulkas2014constructing} apply the bidding strategy of \cite{pandvzic2013offering} for the day-ahead market using stochastic programming for a DH system that includes one CHP unit, a peak boiler and one heat storage tank. 

The so far presented methods focus on the day-ahead market trading only. The consideration of bidding in sequential markets is considered, e.g., in \cite{ayon2017aggregators}, who created bids using stochastic programming in both day-ahead and intra-day markets for an aggregator combining decentralised RES production and consumption without any connection to DH systems. The presented approach first creates bids for the day-ahead market. After this market is cleared, the already committed power production or consumption in the day-ahead market is used to formulate optimal bids for each intra-day market auction throughout the day. Additionally, the standalone participation of different units in sequential electricity markets (especially day-ahead and balancing markets) has been widely discussed in literature (see for instance these sequence bidding strategies for thermal generators \cite{plazas2005multimarket}, microgrids \cite{pei2016optimal}, wind farms \cite{hosseini2013optimal}, hydropower  \cite{vardanyan2017coordinated} or CHP units \cite{kumbartzky2017optimal}). 

To the best of out knowledge, we see a gap regarding the optimal participation of DH systems in a sequential electricity market structure using a realistic framework that includes bidding strategies. There is a need for a planning method that allows DH operators with a portfolio of units to schedule their production under uncertainty and participate in both day-ahead and balancing markets. In particular, for the case when the DH system contains CHP units, power-to-heat technologies and potentially RES power production, which offer the opportunity to lower the heat production costs by trading on the markets. The consideration of all units as a portfolio hedges against the uncertain RES production and resembles the concept of a VPP power producer. However, the operational planning and bidding method needs to account for the limitations of the heat production with respect to demand and thermal storages. Such a method offers the opportunity to analyze the optimal production behaviour of DH systems in a context with RES production. The contributions of this paper can be summarized as follows:
\begin{enumerate}
    \item We bridge the above mentioned gap by extending the VPP bidding method of \cite{pandvzic2013offering} to a DH setting and including balancing market trading explicitly as second step. The underlying stochastic programs are formulated in a general manner to be applicable to arbitrary sets of production units in DH systems.
    \item The method explicitly accounts for the uncertainty coming from RES production in both heat and power and enables us to perform an analysis of the impact of the different uncertainty sources.
    \item We use the method to analyze a real case study based on the Hvide Sande district heating system in Denmark allowing us to draw conclusions on a) the behaviour of the system under uncertain RES production; b) the impact of including balancing market trading to the planning method; c) the benefits of including renewable power production to the portfolio; and d) the annual system costs compared to traditional bidding methods based on forecasts.
    \item An additional contribution is a new approach to generate scenarios for balancing market price scenarios needed for the stochastic programming addressing the balancing market related operation.
\end{enumerate}

Our study is based on the following assumptions. First, we assume the DH operator is a price-taker, i.e., we do not influence the market price, which is reasonable for small- and medium-size DH systems. Second, we assume that the markets allow the submission of price-dependent bids as it is the case in Nordpool. Third, we do not consider minimum and maximum power volume restrictions in both markets. Fourth, we assume electricity prices and RES production are uncertain when planning the day-ahead bidding. For the balancing bids we consider the RES production known for the next hour. The heat demand is assumed to be known and adjusted to cover the heat losses. The reason we consider heat demand as a known parameter is due to its strong correlation to the ambient temperature and season of the year. Thus, having previous observations, we can obtain very accurate predictability 24 hours ahead \cite{petrichenko2017district}. In addition, if a particular deviation from the predicted heat demand occurs, the DH operator have mechanisms to correct these imbalances such as increasing or decreasing the pressure in the DH network. Finally, we do not consider wind spillage as a recourse variable and therefore, we are responsible for our own imbalances.  

The remainder of this paper is organized as follows. In Section \ref{sec:model} we provide the mathematical formulation that describe the two operational problems for day-ahead and balancing market, respectively. Section \ref{sec:uncertainty} describes the modelling of uncertainty, i.e.,  the scenario generation for the RES production and electricity prices. Section \ref{sec:bidding} describes the bidding strategy. The Hvide Sande case study is described in detail in Section \ref{sec:casestudies}. Section \ref{sec:results} provides an analysis and discussion of the results obtained for the case study. Finally, Section \ref{sec:conclusion} summarizes our work and gives an outlook on future work.  

\section{Operational planning model}\label{sec:model}

We start by introducing the two-stage stochastic programs that are the basis for creating bids for the day-ahead and balancing market. The major part of the constraints are valid for both markets and relates to the operation of a portfolio of production units in a district heating system. We start by introducing those constraints. The specific constraints and objectives regarding the two different markets are given in Section \ref{sec:dayaheadmodel} and \ref{sec:balancingmodel} for day-ahead and balancing market, respectively. For an overview of the nomenclature, we refer to Table \ref{tab:nomenclature}. 

\begin{table}
		\footnotesize
	\centering
	\caption{Nomenclature}
	\scalebox{0.97}{
	\begin{tabular}{p{0.17\textwidth}p{0.7\textwidth}}\toprule
		\multicolumn{2}{l}{\textbf{Sets}}\\\midrule
		$\mathcal{T}=\lbrace 1,...,|\mathcal{T}|\rbrace$ & Set of time periods $t$\\
		$\mathcal{U}$ & Set of heat production units $u$\\
		$\mathcal{U}^{\text{CHP}} \subset\mathcal{U} $& Subset of combined heat and power production units\\
		$\mathcal{U}^{\text{H}} \subset\mathcal{U}$ & Subset of heat-only production units \\
		$\mathcal{U}^{\text{EL}} \subset\mathcal{U}$ & Subset of power to heat production units \\ 
		$\mathcal{U}^{\text{RES}} \subset\mathcal{U}$ & Subset of stochastic heat production units \\
		$\mathcal{G}^{\text{}}$ & Set of intermittent renewable power-only producers $g$ \\ 
		$\mathcal{S}^{\text{}}$ & Set of heat storage tanks $s$ \\
		$\Omega$ & Set of scenarios $\omega$\\\midrule
		\multicolumn{2}{l}{\textbf{Parameters}}  \\\midrule
		$C^{\text{H}}_{u}$ & Cost for producing heat with unit $u \in \mathcal{U}$  [DKK/MWh-heat]\\
		$C^{\text{T}}_{g,u}$ & Tariff cost for producing power with unit $g \in \mathcal{G}$ and use it to produce heat in unit $u \in \mathcal{U}^{\text{EL}}$  [DKK/MWh-heat]\\
		$\overline{Q}_{u}/\underline{Q}_{u}$ &Maximum/Minimum heat production for unit $u \in \mathcal{U}$ [MWh-heat]\\
		$A^{\text{DH}}_{u}$ & Binary parameter: 1, if unit $u \in \mathcal{U}$ is connected to the district heating system, 0, otherwise\\
		$A^{\text{S}}_{u,s}$ & Binary parameter: 1, if unit $u \in \mathcal{U}$ is connected to the thermal storage $s$, 0, otherwise\\
		$\varphi_{u}$ &Heat-to-power ratio for unit $u \in \mathcal{U}^{\text{CHP}}$ [$\text{MWh-heat}/\text{MWh-el}$]\\
		$S^{0}_{s}$ & Initial level in storage $s$ [MWh-heat]\\
		$\overline{S^{\text{}}_{s}}/\underline{S^{\text{}}_{s}}$& Maximum/Minimum heat level in storage $s$ [MWh-heat]\\
		${\lambda}_{t}$ &Electricity price for time period $t \in \mathcal{T}$ [DKK/MWh-el]\\
		${\lambda}^{+}_{t}$/${\lambda}^{-}_{t}$ & Penalty for positive/negative imbalance in time period $t \in \mathcal{T}$ [DKK/MWh-el]\\
		${\lambda}^{\text{UP}}_{t}$/${\lambda}^{\text{DOWN}}_{t}$ &Upward/Downward regulating price for time period $t \in \mathcal{T}$ [DKK/MWh-el]\\
		$Q^{\text{D}}_{t}$ &Heat demand for time period $t \in \mathcal{T}$ [MWh-heat]\\
		$P^{\text{RES}}_{g,t,\omega}$ & Stochastic power production of power-only unit $g \in \mathcal{G}^{\text{RES}}$\\
		$Q^{\text{RES}}_{u,t, \omega}$ & Stochastic heat production from heat production unit $u \in \mathcal{U}^{\text{RES}}$\\
		$\pi_{\omega}$ & Probability of scenario $\omega \in \Omega$\\
		$\beta_{}$ & Parameter that determines the deviation of the penalty for the positive and negative imbalance    \\\midrule
		\multicolumn{2}{l}{\textbf{Variables}}  \\\midrule
		$p^{\text{BID}}_{t,\omega} \in \mathbb{R}_{0}$ & Power bid to the day-ahead market unit in period $t \in \mathcal{T}$ [MWh-el] \\
		$q_{u,t,\omega} \in \mathbb{R}^+_{0}$ &Heat production of heat unit $u \in \mathcal{U}$ in period $t \in \mathcal{T}$ [MWh-heat]  \\  
		$q_{u,t,\omega}^{\text{DH}} \in \mathbb {R}^+_{0}$ & Heat production of unit $u \in \mathcal{U}$ inserted to the grid in period $t \in \mathcal{T}$ [MWh-heat] \\  
		$q_{u,s,t,\omega}^{\text{S}} \in \mathbb{R}^+_{0}$ & Heat production of unit $u \in \mathcal{U}$ inserted to the storage $s$ in period $t \in \mathcal{T}$ [MWh-heat]  \\
		$p^{\text{CHP}}_{u,t,\omega} \in \mathbb{R}^+_{0}$ & Power production of unit $u \in \mathcal{U}^{\text{CHP}}$ in period $t \in \mathcal{T}$ [MWh-el] \\
		$p^{\text{GRID}}_{u,t,\omega} \in \mathbb{R}_{0}$ & Power obtained from the grid to produce heat with unit $u \in \mathcal{U}^{\text{EL}}$ in period $t \in \mathcal{T}$ [MWh-el] \\
		$p^{\text{HEAT}}_{g,u,t,\omega} \in \mathbb{R}^+_{0}$ & Power production of unit $g \in \mathcal{G}^{\text{}}$ that serves heat production of unit $u \in \mathcal{U}^{\text{EL}}$ in period $t \in \mathcal{T}$ [MWh-el] \\
		$p^{\text{GEN}}_{g,t,\omega} \in \mathbb{R}^+_{0}$ & Power generation from unit $g \in \mathcal{G}^{\text{}}$ in period $t \in \mathcal{T}$ [MWh-el] \\
	    $p^{\text{+/-}}_{t,\omega} \in \mathbb{R}^+_{0}$ & Positive/Negative power imbalance purchased/sold in period $t \in \mathcal{T}$  and scenario $\omega$ [MWh-el] \\
		$p^{\text{UP/DOWN}}_{t,\omega} \hspace{-0.25cm} \in \mathbb{R}^+_{0}$ & Upward/Downward regulating power purchased/sold in period $t \in \mathcal{T}$  and scenario $\omega$ [MWh-el] \\
		$\sigma_{s,t,\omega} \in \mathbb{R}^+_{0}$ & Level in storage $s$ at time period $t \in \mathcal{T}$ [MWh-heat]  \\
		$\sigma^{\text{OUT}}_{s,t,\omega} \in \mathbb{R}^+_{0}$ & Heat flowing from the storage $s$ to the district heating in period $t \in \mathcal{T}$ [MWh-heat]  \\ \bottomrule
	\end{tabular}}
	\label{tab:nomenclature}
\end{table}

The overall goal is to fulfill the heat demand $Q^{\text{D}}_{t}$ in the district heating network in each period of time $t \in \mathcal{T}$ at lowest cost while taking expected income from bids won on the electricity markets into account. 
The district heating operator has a set of heat and power production units that are operated as portfolio. We divide the set of units in heat producing units $\mathcal{U}$ and intermittent renewable power-only production units $\mathcal{G}$ (wind power or photo-voltaic). The heat producing units $\mathcal{U}$ are further categorized in combined heat and power plants $\mathcal{U}^{\text{CHP}}$ (producing heat and power simultaneously at a heat-to-power ratio $\varphi_{u}$), heat-only units using electricity $\mathcal{U}^{\text{EL}}$, heat-only units with controllable production based on other fuels $\mathcal{U}^{\text{H}}$ and stochastic heat production units $\mathcal{U}^{\text{RES}}$ (e.g. solar thermal).  The stochastic production of both heat and power units are modelled based on a set of scenarios $\Omega$ given by the parameters $Q^{\text{RES}}_{u,t, \omega}$ and $P^{\text{RES}}_{g,t,\omega}$, respectively. Each of the heat producing units has a lower and upper limit on the production amount per period given by $\underline{Q}_{u}$ and $\overline{Q}_{u}$.

The DH operator further uses thermal storages $\mathcal{S}$ to store heat over several periods. The minimum and maximum level of each storage are denoted by $\underline{S^{\text{}}_{s}}$ and $\overline{S^{\text{}}_{s}}$ where as the initial level is set by $S^{0}_{s}$. The physical connections of the units to the storages and the district heating network are modelled by the binary parameters $A^{\text{S}}_{u,s}$ and $A^{\text{DH}}_{u}$, respectively (equals 1, if a connections exists and 0, otherwise).

The operational cost for producing one MWh of heat are represented by the coefficients $C^{\text{H}}_{u}$. A special case is the production of heat based on electricity, i.e., the units $u \in \mathcal{U}^{\text{EL}}$ have additional costs on top of the operational cost based on the electricity needed.  We consider a special tariff $C^{\text{T}}_{g,u}$ for producing heat with heat units $u \in \mathcal{U}^{\text{EL}}$ fueled by power produced by our own power generators $g \in \mathcal{G}$.  Electricity bought from the grid for units $u \in \mathcal{U}^{\text{EL}}$ is included in the bids to the market. The income from the market is approximated based on the amount of power offered to the market and electricity price scenarios $\lambda_{t,\omega}$.

    The decisions determined by the model are the production amounts of heat ($q_{u,t,\omega}$) and power ($p^{\text{CHP}}_{u,t,\omega}$) for the dispatchable units as well as the amount of power offered to the electricity market, the latter being the first-stage decisions in our stochastic program.  Further variables relate to the storage and feeding to the DH and are described later. All variables and their domains are given in Table \ref{tab:nomenclature}.
    
    The following constraints are valid for the production scheduling on both a day-ahead market and balancing market level.
	\begin{subequations}
		{\allowdisplaybreaks
        The heat production of each unit is limited to the capacities of the unit by constraints \eqref{Maxheatproduction}. In constraints \eqref{storagedistrictheating} the production of each unit is split in heat used in the district heating network ($q^{\text{DH}}_{u,t,\omega}$) and heat stored in the thermal storage ($q^{\text{S}}_{u,s,t,\omega}$). The possibility of this split is dependent on the existing connections to storages and the district heating network. Flow in non-existent connections is avoided by constraints \eqref{restrictdh} and \eqref{restricts}.
        
		\begin{align}
		& \underline{Q}_{u,t} \leq q_{u,t,\omega} \leq \overline{Q_{u}} \label{Maxheatproduction} && \forall t \in \mathcal{T}, \forall u \in \mathcal{U}, \forall \omega \in \Omega  \\
		&  q_{u,t, \omega} = q^{\text{DH}}_{u,t,\omega}+\sum_{s \in \mathcal{S}}q^{\text{S}}_{u,s,t,\omega}  \label{storagedistrictheating} && \forall t \in \mathcal{T}, \forall u \in \mathcal{U}, \forall \omega \in \Omega\\
		&  q^{\text{DH}}_{u,t,\omega} \leq \overline{Q_{u}}A^{\text{DH}}_{u}  \label{restrictdh} && \forall t \in \mathcal{T}, \forall u \in \mathcal{U},  \forall \omega \in \Omega\\
		&  q^{\text{S}}_{u,s,t,\omega}\leq \overline{Q_{u}}A^{\text{S}}_{u,s} \label{restricts} && \forall t \in \mathcal{T}, \forall u \in \mathcal{U}, \forall \omega \in \Omega
		\end{align}
        The coupling of heat and power production in CHP units is modelled in constraints \eqref{Heattopowerchp}. Furthermore, the electric boiler production can be based on electricity bought on the market ($p^\text{GRID}_{u,t,\omega}$) or from our own power generators ($p^\text{HEAT}_{g,u,t,\omega}$) (see constraints \eqref{Heattopowerel}). Stochastic renewable heat production from, e.g, solar thermal units, is dependent on the scenario and given as input in constraints \eqref{resheat}.
        
		\begin{align}
		& q_{u,t,\omega} = \varphi_{u}p^{\text{CHP}}_{u,t\omega} \label{Heattopowerchp}  &&  \forall t \in \mathcal{T}, \forall u \in \mathcal{U}^{\text{CHP}}, \forall \omega \in \Omega  \\
			& q_{u,t,\omega} = \varphi_{u}\big(p^{\text{GRID}}_{u,t,\omega}+\sum_{g \in \mathcal{G}}p^{\text{HEAT}}_{g,u,t,\omega}\big) \label{Heattopowerel}  &&  \forall t \in \mathcal{T},  \forall u \in \mathcal{U}^{\text{EL}}, \forall \omega \in \Omega \\
			& q_{u,t,\omega} = Q^{\text{RES}}_{u,t,\omega}  &&  \forall t \in \mathcal{T}, \forall u \in \mathcal{U}^{\text{RES}}, \forall \omega \in \Omega  \label{resheat}
		\end{align}
		The thermal storage level ($\sigma_{s,t,\omega}$) limitations as well as in- and outflows ($\sigma^{OUT}_{s,t,\omega}$) are modelled in constraints \eqref{eq2:thermal_limits} and \eqref{eq2:thermal_level}, respectively. At the end of the planning horizon, we impose that the storage level is at least as high as in the beginning of the planning horizon to avoid emptying the storage in every optimization (constraints \eqref{eq2:thermal_end}).
		
		\begin{align}
		&\underline{S}_{s}  \leq \sigma_{s,t,\omega} \leq \overline{S}_{s}  && \forall t \in \mathcal{T}, \forall s \in \mathcal{S}, \forall \omega \in \Omega \label{eq2:thermal_limits}\\
			& \sigma_{s,t,\omega} =\sigma_{s,t-1,\omega}+\sum_{u \in \mathcal{U}_{}}q^{\text{S}}_{u,s,t,\omega}-\sigma^{\text{OUT}}_{s,t,\omega}  && \forall t \in \mathcal{T}, \forall s \in \mathcal{S}, \forall \omega \in \Omega
		\label{eq2:thermal_level}\\
		&\sigma_{s,|\mathcal{T}|,\omega}   \geq S^{0}_{s}  && \forall s \in \mathcal{S}, \forall \omega \in \Omega \label{eq2:thermal_end} 
		\end{align}
        The heat demand in the network in each period is ensured by constraints \eqref{eq2:Heatbalance} by using either heat directly from the units or from the storage.
        
		\begin{align}
		& Q^{\text{D}}_{t} = \sum_{u \in \mathcal{U}}q^{\text{DH}}_{u,t,\omega}+\sum_{s \in \mathcal{S}}s^{\text{OUT}}_{s,t,\omega}  && \forall t \in \mathcal{T}, \forall \omega \in \Omega \label{eq2:Heatbalance}
		\end{align}
		The renewable power production from the stochastic power generators is modelled in constraints \eqref{respower} depending on the scenario. The power can be used either to produce heat with the electric boiler ($p^{\text{HEAT}}_{g,u,t,\omega}$) or sold on the market ($p^{\text{GEN}}_{g,t,\omega}$).
		
		\begin{align}
		&p^{\text{GEN}}_{g,t,\omega}+\sum_{u \in \mathcal{U}^{\text{EL}}}p^{\text{HEAT}}_{g,u,t,\omega} =  P^{\text{RES}}_{g,t,\omega}  &&  
		 \forall t \in \mathcal{T}, \forall g \in \mathcal{G}^{},  \forall \omega \in \Omega\label{respower}
		\end{align}}
	\end{subequations}
Based on this initial set of constraints, the model is extended for day-ahead or balancing market optimization in the succeeding sections.

\subsection{Optimization for the day-ahead market}\label{sec:dayaheadmodel}

The first-stage variables (here-and-now decisions) for the day-ahead market production scheduling are the power bids $p^{\text{BID}}_{t,\omega}$ for each hour of the next day $t \in \lbrace 1,\dots, 24\rbrace$. As these are dependent on the production of all other dispatchable units, we determine the heat ($q_{u,t,\omega}$) and power production ($p_{u,t,\omega}$) of all units  as well as the power bid amounts for the remaining planning horizon ($p^{\text{BID}}_{u,t,\omega} \forall t \in \lbrace 25,\dots,|\mathcal{T}|\rbrace$) as second stage variables.

	\begin{subequations}
The objective function \eqref{objectivefunction1} minimizes the expected cost of producing heat by all units minus the expected income for the day-ahead electricity market. Deviations from the day-ahead market bid are penalized by paying the imbalances ($p^{+}_{t,\omega}, p^{-}_{t,\omega}$) at the balancing stage.

		\begin{align}
		\underset{}{\text{min}} & \hspace{0.2cm} \sum_{t \in \mathcal{T}}\sum_{\omega \in \Omega} \pi_{\omega} \Bigg[  \sum_{u \in \mathcal{U}^{\text{CHP}}}C^{\text{H}}_{u}q_{u,t,\omega}+ \sum_{u \in \mathcal{U}^{\text{H}}} C_{u}^{\text{H}}q_{u,t,\omega} + \sum_{u \in \mathcal{U}^{\text{EL}}}C_{u}^{\text{H}}p^{\text{GRID}}_{u,t,\omega}\notag  \\ 
		 & +\sum_{g \in \mathcal{G}} \sum_{u \in \mathcal{U}^{\text{EL}}}C^{T}_{g,u}p^{\text{HEAT}}_{g,u,t,\omega} 
		   -\big(\lambda_{t,\omega}p^{\text{BID}}_{t,\omega} - \lambda^{+}_{t}p^{+}_{t,\omega}+\lambda^{-}_{t}p^{-}_{t,\omega}\big) \Bigg]\label{objectivefunction1}
		\end{align}
The bidding amount $p^{\text{BID}}_{t,\omega}$ is dependent on the power production from CHP units and the generator as well as the power used for the electric boiler (see constraints \eqref{totalpowersold}). Any deviations from the bidding amount are captured in the variables $p^{+}_{t,\omega}$ and $p^{-}_{t,\omega}$ to be penalized in the objective function.

	\begin{align}
		& p^{\text{BID}}_{t,\omega} = \negthickspace \negthickspace \negthickspace \sum_{u \in \mathcal{U}^{\text{CHP}}} \negthickspace \negthickspace p^{\text{CHP}}_{u,t,\omega}+\sum_{g \in \mathcal{G}^{\text{}}}p^{\text{GEN}}_{g,t,\omega}-\negthickspace \negthickspace \negthickspace\sum_{u \in \mathcal{U}^{\text{EL}}}\negthickspace \negthickspace p^{\text{GRID}}_{u,t,\omega}+p^{+}_{t,\omega}-p^{-}_{t,\omega} &&  
		 \forall t \in \mathcal{T}, \forall \omega \in \Omega \label{totalpowersold}
	\end{align}
	The equations \eqref{equalitiybids} are based on the method in \cite{pandvzic2013offering} and ensure that only one bidding curve, i.e., one set of power amount and price pairs, is created per time period $t$ while constraints \eqref{nondecreasingbids} ensure that the bidding curves are non-decreasing for all time steps $t \in \mathcal{T}$. 	
	
				\begin{align}
	    & p^{\text{BID}}_{t,\omega} = p^{\text{BID}}_{t,\omega'} && \forall t \in \mathcal{T}, \forall (\omega, \omega') \in \Omega && : \lambda_{t,\omega}=\lambda_{t,\omega'} \label{equalitiybids} \\
	     & p^{\text{BID}}_{t,\omega} \leq p^{\text{BID}}_{t,\omega'} && \forall t \in \mathcal{T}, \forall (\omega, \omega') \in \Omega && : \lambda_{t,\omega} \leq \lambda_{t,\omega'} \label{nondecreasingbids}
		\end{align}
	\end{subequations}
	 The operational model to optimize the production for the day-ahead market bidding can be summarized as follows in \eqref{eq:dayahead1} to \eqref{eq:dayahead3}. 
	 
		\begin{subequations}
			\begin{align}
	     \text{min} \hspace{0.2cm} &  \eqref{objectivefunction1} \label{eq:dayahead1}\\
	     \text{s.t}  \hspace{0.2cm}&  \eqref{Maxheatproduction}-\eqref{respower} \label{eq:dayahead2}\\
	                \hspace{0.2cm}  & \eqref{totalpowersold}-\eqref{nondecreasingbids}\label{eq:dayahead3}
		\end{align}
	\end{subequations}
		To avoid speculation in the operation of the system, we define the penalty costs for deviation as follows. 
\begin{gather*}
\lambda^{+}_{t,\omega}=
\begin{cases}
  \lambda^{}_{t,\omega}+\beta \cdot \lambda^{}_{t,\omega} & \text{if } \lambda^{}_{t,\omega}\geq 0 \\    
  \lambda^{}_{t,\omega}-\beta \cdot \lambda^{}_{t,\omega} & \text{if } \lambda^{}_{t,\omega}< 0 \\     
\end{cases};	
\hspace{0.3cm}
\lambda^{-}_{t,\omega}=
\begin{cases}
  \lambda^{}_{t,\omega}-\beta \cdot \lambda^{}_{t,\omega} & \text{if } \lambda^{}_{t,\omega}\geq 0 \\    
  \lambda^{}_{t,\omega}+\beta \cdot \lambda^{}_{t,\omega} & \text{if } \lambda^{}_{t,\omega}< 0 \\     
\end{cases}	
\end{gather*}
where $\beta$ is a parameter with value greater than 0. Thus, we ensure that the penalty to pay would be higher than the day-ahead prices in case of positive deviation. On the contrary, in case of producing more power than sold in the day-ahead market, the profits for selling that excess power on the balancing market are always lower than selling that energy in the day-ahead market. Therefore, the model tries to sell the right amount of power on the day-ahead market and avoid imbalances.    

\subsection{Optimization for the balancing market}\label{sec:balancingmodel}

The balancing market problem is solved once per hour and like in the day-ahead problem \eqref{eq:dayahead1}-\eqref{eq:dayahead3}, we generate non-decreasing bidding curves using the stochastic formulation of the problem. In this case, the first-stage decisions are the upward ($p^{\text{UP}}_{t,\omega}$) and downward ($p^{\text{DOWN}}_{t,\omega}$) regulation offered to formulate the bidding curves for the balancing market. The remaining variables can be adapted to the realization of the uncertainty and considered as second-stage decisions. In this formulation of the balancing problem, the committed power production or consumption for the day-ahead is given as a parameter ($\widehat{p}^{\text{BID}}_{t,\omega}$). Due to the high unpredictability of the balancing prices we use $\mathcal{T^{\text{B}}}$ periods as the planning horizon for the balancing problem, which can be shorter than the horizon used in the day-ahead problem. Upward regulation ($p^{\text{UP}}_{t,\omega}$) is provided in case there is a need for more power in the system, therefore the producer has the opportunity to sell additional power at the upward regulating price ($\lambda^{\text{UP}}_{t,\omega}$). On the contrary, if the systems has excess of production, the TSO activates offers for downward regulation, where producers can consume power ($p^{\text{DOWN}}_{t,\omega}$) at the downward regulating price ($\lambda^{\text{DOWN}}_{t,\omega}$).

The objective function \eqref{objectivefunction2} for the balancing problem again minimizes the cost considering income from the market and penalties for imbalances. 

\begin{subequations}
\begin{align}
    		\underset{}{\text{min}} & \hspace{0.2cm}  \sum_{t \in \mathcal{T^{\text{B}}}}\sum_{\omega \in \Omega} \pi_{\omega} \Bigg[ \sum_{u \in \mathcal{U}^{\text{CHP}}}C^{\text{H}}_{u}q_{u,t,\omega}+  \sum_{u \in \mathcal{U}^{\text{H}}} C_{u}^{\text{H}}q_{u,t,\omega} + \sum_{u \in \mathcal{U}^{\text{EL}}}C_{u}^{\text{H}}p^{\text{GRID}}_{u,t,\omega}  \label{objectivefunction2}\\ 
		 & +\sum_{g \in \mathcal{G}} \sum_{u \in \mathcal{U}^{\text{EL}}}C^{T}_{g,u}p^{\text{HEAT}}_{g,u,t,\omega} 
		   - \big( \lambda^{-}_{t}p^{-}_{t,\omega}- \lambda^{+}_{t}p^{+}_{t,\omega}+\lambda^{\text{UP}}_{t,\omega}p^{\text{UP}}_{t,\omega}-\lambda^{\text{DOWN}}_{t,\omega}p^{\text{DOWN}}_{t,\omega}\big) \Bigg]\notag
\end{align}
 The balance in the power production is ensured in equations \eqref{eq:powerbalance2}. Here the power committed on the day-ahead market is given as a parameter ($\widehat{p}^{\text{BID}}_{t,\omega}$). To balance the production with the bidding amount, constraint \eqref{eq:powerbalance2} can either use the variables determining the upward ($p^{\text{UP}}_{t,\omega}$) or downward regulation ($p^{\text{DOWN}}_{t,\omega}$) amounts or pay imbalances. The imbalances are captured in  $p^{+}_{t,\omega}$ and $p^{-}_{t,\omega}$.
 
    	\begin{align}
		& \widehat{p}^{\text{BID}}_{t,\omega} = \negthickspace \negthickspace \negthickspace \sum_{u \in \mathcal{U}^{\text{CHP}}} \negthickspace \negthickspace p^{\text{CHP}}_{u,t,\omega}+\sum_{g \in \mathcal{G}^{\text{}}}p^{\text{GEN}}_{g,t,\omega}-\negthickspace \negthickspace \negthickspace\sum_{u \in \mathcal{U}^{\text{EL}}}\negthickspace \negthickspace p^{\text{GRID}}_{u,t,\omega}+p^{+}_{t,\omega}-p^{-}_{t,\omega}-p^{\text{UP}}_{t,\omega}+p^{\text{DOWN}}_{t,\omega} \label{eq:powerbalance2} \\  
		&   \forall t \in \mathcal{T^{\text{B}}}, \forall \omega \in \Omega \notag
\end{align}
To ensure ordered bidding curves in the balancing market, we define constraints \eqref{equalitiybidsbal1} and \eqref{equalitiybidsbal2} analogously to the day-ahead market problem. Here the offers for upward regulation and downward regulation, present a non-decreasing and non-increasing order, respectively.

\begin{align}
 & p^{\text{UP}}_{t,\omega} \leq p^{\text{UP}}_{t,\omega'} && \forall t \in \mathcal{T^{\text{B}}}, \forall (\omega, \omega') \in \Omega && : \lambda^{\text{UP}}_{t,\omega} \leq \lambda^{\text{UP}}_{t,\omega'} \label{equalitiybidsbal1} \\
 & p^{\text{DOWN}}_{t,\omega} \geq p^{\text{DOWN}}_{t,\omega'} && \forall t \in \mathcal{T^{\text{B}}}, \forall (\omega, \omega') \in \Omega && : \lambda^{\text{DOWN}}_{t,\omega} \leq \lambda^{\text{DOWN}}_{t,\omega'} \label{equalitiybidsbal2}
\end{align}
\end{subequations}
The entire formulation for the balancing market problem is given by \eqref{eq:balancing1} to \eqref{eq:balancing3}.

\begin{subequations}
\begin{align}
	     \text{min} \hspace{0.2cm} &  \eqref{eq:powerbalance2} \label{eq:balancing1}\\
	     \text{s.t}  \hspace{0.2cm}&  \eqref{Maxheatproduction}-\eqref{respower} \label{eq:balancing2}\\
	                \hspace{0.2cm}  & \eqref{eq:powerbalance2}-\eqref{equalitiybidsbal2}\label{eq:balancing3}
\end{align}
Furthermore, as in the day-ahead problem, we need to prohibit speculation of  the system by defining the penalty prices $\lambda^{+}_{t,\omega}$ and $\lambda^{-}_{t,\omega}$ as follows. 

\begin{gather*}
\lambda^{+}_{t,\omega}=
\begin{cases}
  \lambda^{}_{t,\omega}+\beta \cdot \lambda^{}_{t,\omega} & \text{if } \lambda^{}_{t,\omega}\geq 0, p^{\text{UP}}_{t,\omega} = 0\\    
  \lambda^{}_{t,\omega}-\beta \cdot \lambda^{}_{t,\omega} & \text{if } \lambda^{}_{t,\omega}< 0 , p^{\text{UP}}_{t,\omega} = 0\\ 
  \lambda^{\text{UP}}_{t,\omega}+\beta \cdot \lambda^{\text{UP}}_{t,\omega} & \text{if } \lambda^{}_{t,\omega}\geq 0, p^{\text{UP}}_{t,\omega} \geq 0 \\    
  \lambda^{\text{UP}}_{t,\omega}-\beta \cdot \lambda^{\text{UP}}_{t,\omega} & \text{if } \lambda^{}_{t,\omega}< 0, p^{\text{UP}}_{t,\omega} \geq 0\\
  \end{cases}\\
\lambda^{-}_{t,\omega}=
\begin{cases}
  \lambda^{}_{t,\omega}-\beta \cdot \lambda^{}_{t,\omega} & \text{if } \lambda^{}_{t,\omega}\geq 0, p^{\text{DOWN}}_{t,\omega} = 0\\    
  \lambda^{}_{t,\omega}+\beta \cdot \lambda^{}_{t,\omega} & \text{if } \lambda^{}_{t,\omega}< 0, p^{\text{DOWN}}_{t,\omega} = 0\\ 
  \lambda^{\text{DOWN}}_{t,\omega}-\beta \cdot \lambda^{\text{DOWN}}_{t,\omega} & \text{if } \lambda^{}_{t,\omega}\geq 0, p^{\text{DOWN}}_{t,\omega} \geq 0 \\    
  \lambda^{\text{DOWN}}_{t,\omega}+\beta \cdot \lambda^{\text{DOWN}}_{t,\omega} & \text{if } \lambda^{}_{t,\omega}< 0, p^{\text{DOWN}}_{t,\omega} \geq 0\\    
\end{cases}	
\end{gather*}
\end{subequations}

\section{Modeling Uncertainty}\label{sec:uncertainty}

In particular the day-ahead market optimization includes uncertainty with respect to the production of the stochastic production units (wind power and solar thermal). But both planning problems also have to consider that the electricity prices are still unknown at the time of planning. To account for these uncertainties, we include them as scenarios to our two-stage stochastic programs. The remainder of this section describes the forecasting and scenario generation process. 

\subsection{Wind power production forecast}\label{sec:scenarios_wind}

\begin{figure}
	\begin{subfigure}{.48\textwidth}
	\centering
	\includegraphics[width=0.8\textwidth]{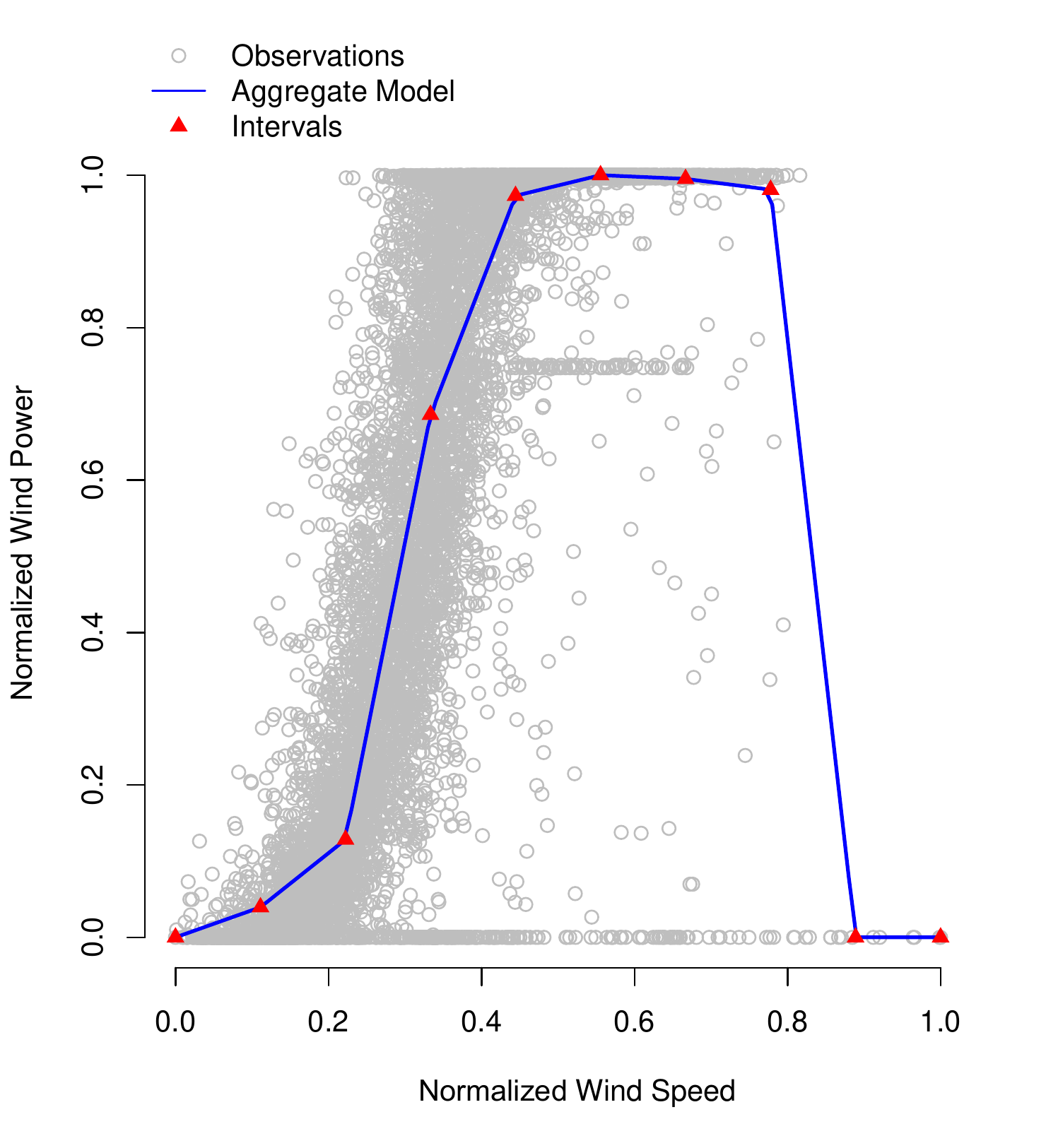}
	\caption{Wind power curve using real data for one year power production and wind speed in NordPool DK1}
	\label{WindPowerCurve}
	\end{subfigure}
	\quad
	\begin{subfigure}{.48\textwidth}
	\centering
	\includegraphics[width=0.8\textwidth]{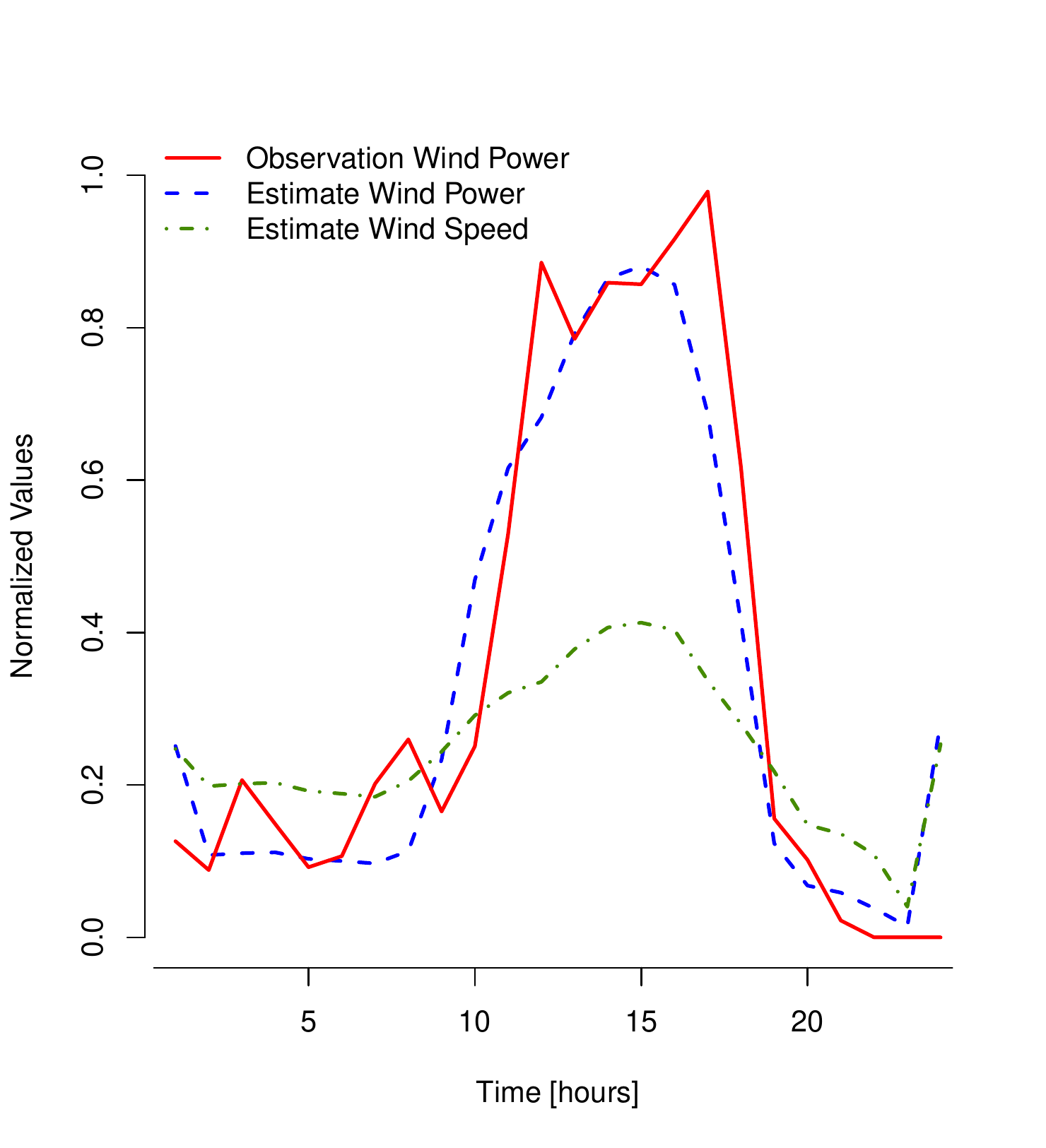}
	\caption{Wind power predictions for one day receding horizon and the normalized wind speed}
	\label{WindPowerPredictions}
	\end{subfigure}
	\caption{Wind power prediction process} \label{WindForecast}
\end{figure}

For an easy replicability of our experiments, we use a wind forecast based on local linear regressions of the wind power curve \cite{pinson2008local}. As Figure \ref{WindPowerCurve} shows, the power curve is divided into intervals with equal distribution based on the normalized wind speed. For each interval, a linear regression is fitted to the data using a least squares estimate. The linear regressions are later integrated into one single function. From this aggregated function, we can predict the wind power production using the wind speed forecast as depicted in Figure \ref{WindPowerPredictions}.

\subsection{Solar Thermal Forecast}\label{sec:scenarios_solar}

The appropriate function to predict solar thermal forecast depends on the technology used in the solar collectors. In this work, we consider flat thermal solar collectors with a fixed inclination angle and orientated towards maximizing the solar radiation during the summer season. The forecasting technique used here is presented in \cite{Solarcol35:online} and given in \eqref{eq:solarpred}.
\begin{align}
    Q_{t} = A^{\text{S}}\Big[I^{\text{D}}_{t}\gamma-\eta_{1}\big(T^{\text{AVG}}_{t}-T^{\text{AMB}}_{t}\big)-\eta_{2}\big(T^{\text{AVG}}_{t}-T^{\text{AMB}}_{t}\big)^{2}\Big] & \hspace{0.1cm} \forall t \in \mathcal{T}\label{eq:solarpred}
\end{align}
where $Q_{t}$ is the heat production at time $t$, $A^{\text{S}}$ is the area of the entire solar thermal field and $I_{t}^{\text{D}}$ is the solar radiation (including direct and diffusive) that heats the solar collectors for time period $t$. $T_{t}^{\text{AVG}}$ and $T_{t}^{\text{AMB}}$ are the average temperature inside the solar collector and the outside temperature, respectively. The remaining parameters ($\gamma$,$\eta_{1}$,$\eta_{2}$) are the coefficients of the equations. The average temperature ($T_{t}^{\text{AVG}}$) is defined as the average between the cold water entering and the hot water leaving the solar collector. For the sake of simplicity, we consider this temperature as constant $\forall t \in \mathcal{T}$.

\subsection{Day-ahead electricity price forecast}\label{sec:scenarios_dayahead}

Electricity prices in day-ahead markets present an autocorrelation and seasonal variation that usually can be detect using time series models. For this work, the electricity price forecast is obtained using a SARMAX model with a daily seasonality pattern that has been successfully applied to predict electricity prices \cite{gonzalez2012forecasting}. In addition, an exogenous variable based on Fourier series is used to describe the weekly seasonality \cite{ringwood1993forecasting}. This results in the following model \eqref{eq:sarmax}.

\begin{subequations}
\begin{align}
     \lambda_{t} = \mu + \phi_{1} \lambda_{t-1}+\phi_{2}\lambda_{t-2}+\phi_{24}\lambda_{t-24}+\theta_{1}\varepsilon_{t-1}+\theta_{2}\varepsilon_{t-2}+\theta_{24}\varepsilon_{t-24}+X   \label{eq:sarmax}
\end{align}
The estimated electricity price ($\lambda_{t}$) for time period $t$ is calculated by the linear combination of the intercept $\mu$, the autoregressive (AR) terms $\lambda_{t-1}$, $\lambda_{t-2}$ and $\lambda_{t-24}$ and the moving average (MA) terms $\varepsilon_{t-1}$, $\varepsilon_{t-2}$ and $\varepsilon_{t-24}$ for 1, 2 and 24 hours prior to time period $t$. The forecast parameters $\phi_{1}$, $\phi_{2}$, $\phi_{24}$, $\theta_{1}$, $\theta_{2}$ and $\theta_{24}$ are updated on a daily basis. The exogenous variable $X$ allows to integrate external variables into the model, in our case the Fourier series describing the weekly seasonality of the data \eqref{eq:fourier}. 

\begin{align}
    X=\sum_{k=1}^{K} \alpha_{k} \sin{\bigg(\frac{2\pi kt}{T}\bigg)}+\sum_{k=1}^{K} \beta_{k} \cos{\bigg(\frac{2\pi kt}{T}\bigg)}\label{eq:fourier}
\end{align}
\end{subequations}
where $K$ determines the number of Fourier terms considered (chosen by minimizing the AICc value). The parameter $T$ represents the seasonality period in the series, in our case we consider a weekly seasonality of $T=168$. Finally, $\alpha_{t}$ and $\beta_{t}$ represent the forecast parameters for the weekly seasonality, and like the forecast parameters for the AR and MA terms, both are updated on a daily basis. 

\subsection{Scenario generation for RES production and day-ahead market prices}\label{sec:scenario_generation}

The forecasts for the three previously mentioned data sets are based on probabilistic forecasts. Therefore, we generate scenarios using a Monte-Carlo simulation applying a multivariate Gaussian distribution with zero mean that describes the stochastic process, which we consider as stationary, in our predictions. We use the algorithm presented in \cite{conejo2010decision} to initialize the scenario generation process and randomly generate the error terms. The algorithm is repeated for each time period in the receding horizon and for all scenarios. In our case, we generate a random walk for the time horizon using normalized white noise that we iteratively add to the predicted value resulting in one scenario. 

To get a representative set of scenarios, we generate a large amount of equiprobable scenarios. Those are reduced to the desired number by applying the clustering technique \textit{partition around medoids} (PAM) \cite{reynolds2006clustering}. Each medoid scenario is a scenario in our model, while the probability is obtained by the sum of the scenarios attached to the medoid.


\subsection{Scenario generation for balancing prices}\label{sec:scenarios_balancing}

The generation of scenarios for balancing prices is less intuitive compared to the day-ahead market prices described before. In particular, because there is not always a need for upward or downward regulation, and if there is, the regulating prices are defined as a function of the imbalanced power volume which makes these prices very hard to predict. The method proposed in \cite{olsson2008modeling} is widely used in literature to create balancing price forecasts. The authors develop a model that combines a SARIMA to predict the amount of upward and downward regulating prices in combination with a discrete Markov model representing the discontinuous variability in the activation of upward and downward regulation. This variability is represented through a matrix that indicates the transition probability between states. Using this techniques, scenarios can be generated by sampling the error term in the time series models and creating different sequences for the Markov model. 

In this section, we propose a novel approach to generate balancing prices scenarios. Our motivation to use a different new scenario generation technique for real-time balancing prices is due to the fact that the authors in \cite{olsson2008modeling} apply their method in a specific bidding area where prices follow a regular shape and pattern that can be accurately predicted, i.e., regions with low integration of RES. In systems with a high penetration of RES (especially wind power), large imbalances can occur in a very short time and thereby affect the balancing prices, which respond to the volume of the imbalance. Due to this variability, balancing prices do not necessarily follow a trend that can be easily predicted using time series models. Furthermore, the method proposed by \cite{olsson2008modeling} models the probability of imbalance states and does not consider the specific duration of these states. We think that this duration must be taken into account since the upward and downward regulation prices are affected by this duration. 
 
Our approach is based on the algorithm to create unit availability scenarios presented in \cite{conejo2010decision}. Initially, the following methodology is applied for upward and downward regulation separately. The results are combined in a final step. The generation of the final predicted prices is carried out based on sampling the deviation compared to the day-ahead price (in \%). 

\begin{figure}
	\begin{subfigure}{.49\textwidth}
	\centering
	\includegraphics[width=0.8\textwidth]{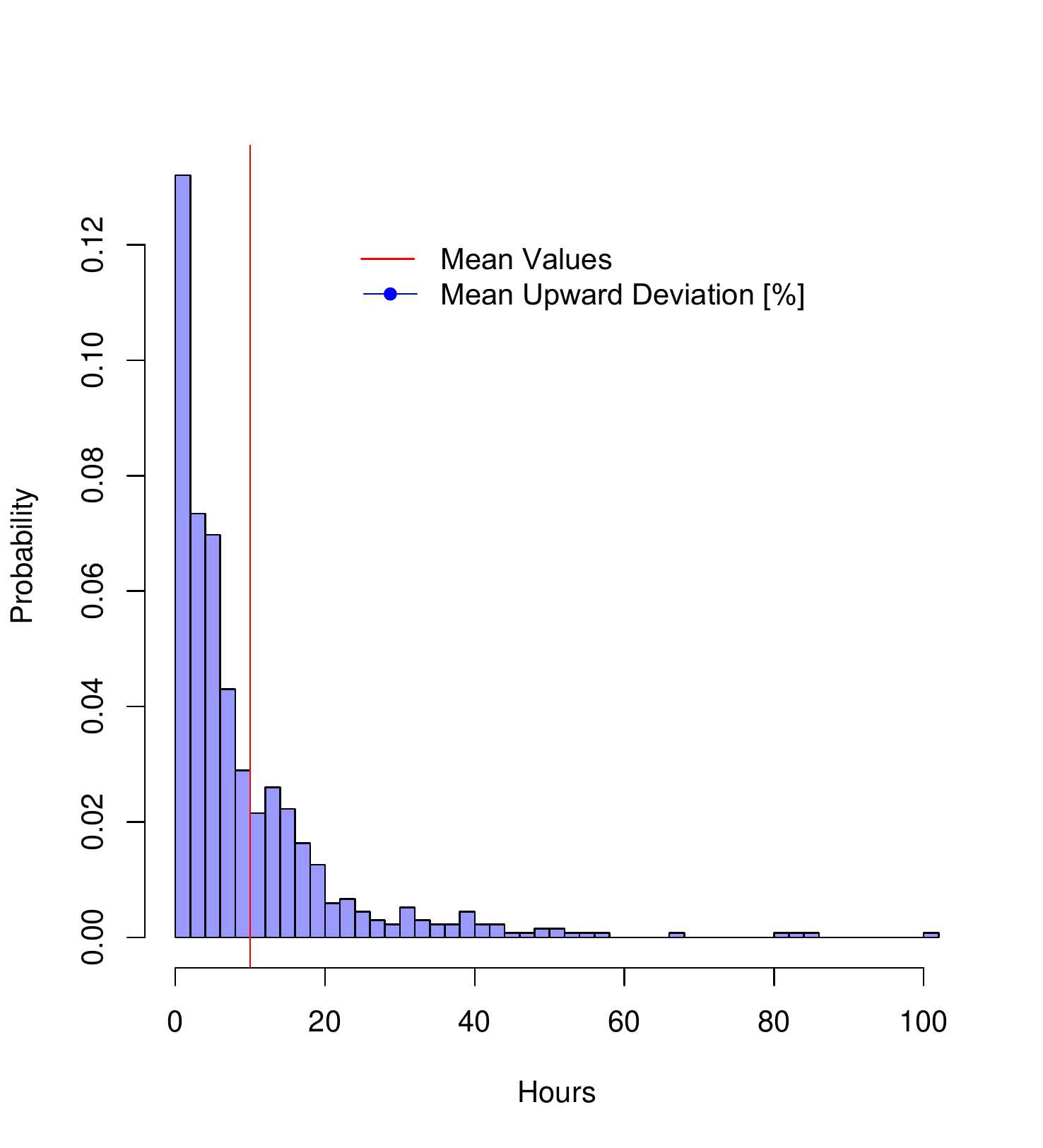}
	\caption{Time elapsed between upward regulation periods}
	\label{fig:timetoactiveupward}
	\end{subfigure}
	\begin{subfigure}{.49\textwidth}
	\centering
	\includegraphics[width=0.8\textwidth]{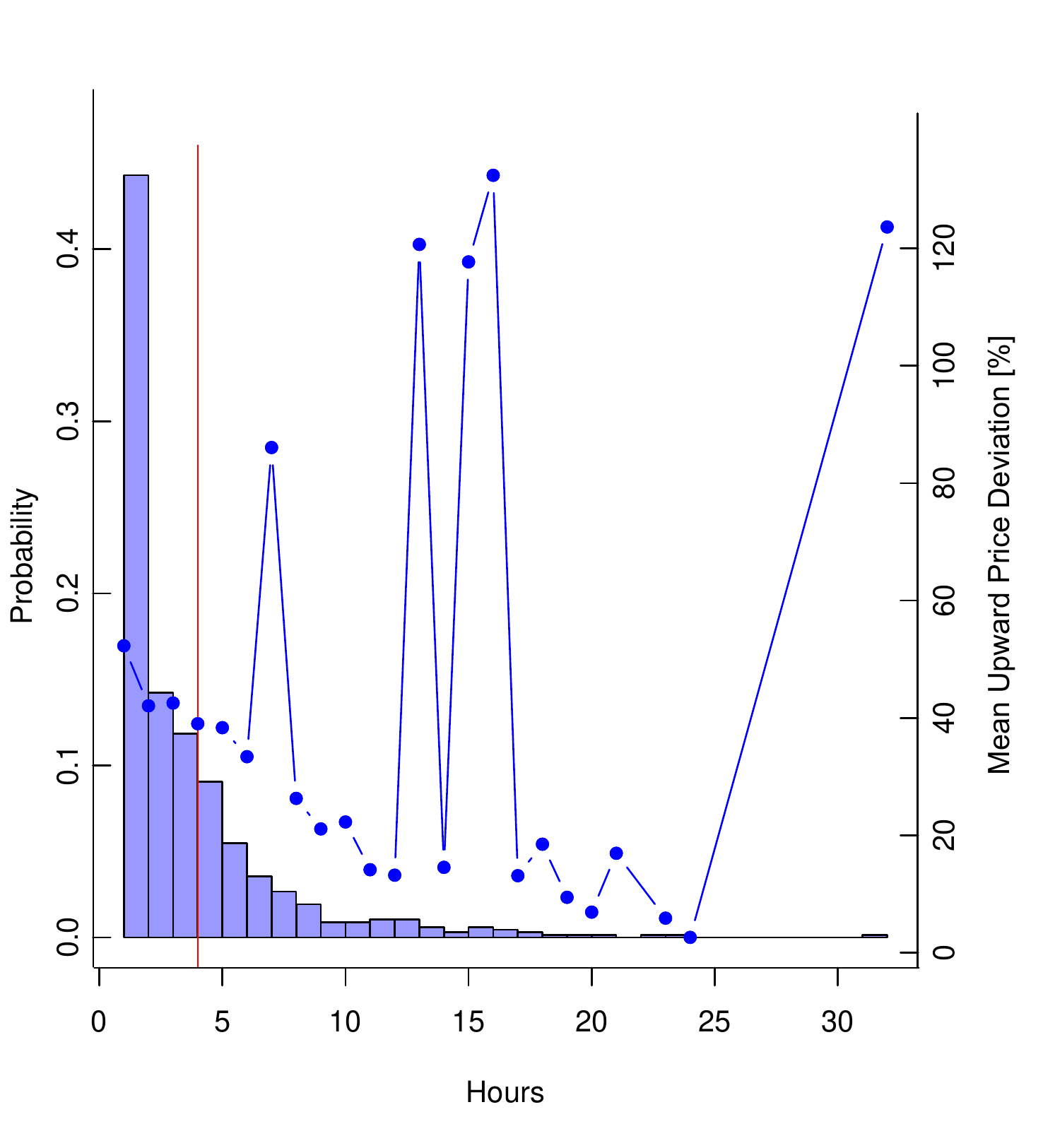}
	\caption{Duration of upward regulation and price deviation}
	\label{fig:DurationUpward}
	\end{subfigure}
	\caption{Distributions of elapsed time between and duration of upward regulation as well as average regulating prices for year 2017 in the NordPool bidding area DK1} \label{fig:Scenariogenerationbalancing}
\end{figure}

The first step is to gather previous observations from the balancing market to determine the experimental distribution of the duration (time elapsed) in between two upward regulation periods or downward regulation periods, respectively, and the corresponding mean values $\tau^{\text{T}+}$ and $\tau^{\text{T}-}$. An example for upward regulation is given in Figure \ref{fig:timetoactiveupward}, where the red line represents the mean value. 
In addition, the distribution of the actual duration for each upward and downward regulation period is also obtained (see Figure \ref{fig:DurationUpward} for upward regulation) along with the mean duration $\tau^{\text{D}+}$ and $\tau^{\text{D}-}$. At the same time, the observed deviations between day-ahead and balancing market prices are averaged for each duration of regulation (see function in Figure \ref{fig:DurationUpward}). By connecting those mean duration values, we get the functions $f^{+}(x)$ and $f^{-}(x)$ telling us for each duration of regulation the  deviation from the day-ahead market price for upward and downward regulation prices, respectively.

Once the experimental distribution and values for $\tau^{\text{T}+}$, $\tau^{\text{T}-}$, $\tau^{\text{D}+}$, $\tau^{\text{D}-}$, $f^{+}(\tau^{\text{D}+})$ and $f^{-}(\tau^{\text{D}-})$ are obtained, the scenario generation is started. As in \cite{conejo2010decision}, we assume that $\tau^{\text{T}+}$, $\tau^{\text{T}-}$, $\tau^{\text{D}+}$ and $\tau^{\text{D}-}$ can be characterized as random variables that follow an exponential distribution, which is a reasonable assumption confirmed by the observations shown in Figure \ref{fig:Scenariogenerationbalancing}. Therefore, random samples of these values can be obtained by applying equations \eqref{randomsample2}, where $u_{1}$ and $u_{2}$ are uniformly distributed variables between 0 and 1.

\begin{equation}
    \tau^{\text{T}(+/-)}_{\omega}=-\tau^{\text{T}(+/-)}\cdot\ln(u_{1}); \hspace{0.5cm}
    \tau^{\text{D}(+/-)}_{\omega}=-\tau^{\text{D}(+/-)}\cdot\ln(u_{2}) \label{randomsample2}
\end{equation}

\begin{algorithm}[t]
\caption{Generate balancing price scenarios}
		\footnotesize
	\begin{algorithmic}[1]
		\For{each $\omega \in \Omega$}
		\State $t \gets 1$
		\While{$t \leq |\mathcal{T}|$}
		\State $\tau^{\text{T}(+/-)}_{\omega}$=-$\tau^{\text{T}(+/-)}\cdot\ln(u_{1})$ where $u_{1} \sim \mathcal{U}(0,1)$ is random
		\State $\tau^{\text{D}(+/-)}_{\omega}$=-$\tau^{\text{D}(+/-)}\cdot\ln(u_{2})$ where $u_{2}\sim \mathcal{U}(0,1)$ is random
		\State $t^{\text{Start}} \gets \text{min}\{|\mathcal{T}|,round(t+\tau^{\text{T}(+/-)}_{\omega})\} $ 
		\State $t^{\text{End}} \gets \text{min}\{|\mathcal{T}|,round(t+\tau^{\text{T}(+/-)}_{\omega} + \tau^{\text{D}(+/-)}_{\omega})\} $ 
		\For{$t' = t$ to $t^{\text{Start}}$}
		\State $\Delta \lambda^{(\text{UP}/\text{DOWN})}_{t',\omega}$ = 0
		\EndFor
		\For{$t' = t^{\text{Start}}+1$ to $t^{\text{End}} $}
		\State $\Delta \lambda^{(\text{UP}/\text{DOWN})}_{t',\omega}$ = $f^{(+/-)}(\tau^{\text{D}(+/-)})+\varepsilon^{(+/-)}_{t'}$ where $\varepsilon^{(+/-)}_{t'} \sim \mathcal{N}(\mu,\sigma^{2})$ is random
		\EndFor
		\State $t \gets t^{\text{End}} + 1$
		\EndWhile
		\EndFor
		\State Return $\Delta \lambda^{(+/-)}_{t,\omega}$
	\end{algorithmic}
	\label{alg:balancingprices}
\end{algorithm}

The algorithm to generate $|\Omega|$ with a time horizon of $|\mathcal{T}|$ periods is summarized in Algorithm \ref{alg:balancingprices} and works as follows. For each scenario we move through the forecasting horizon starting at period 1. The time to the next regulation period $\tau^{\text{T}(+/-)}_{\omega}$ and the duration of this period $\tau^{\text{D}(+/-)}_{\omega}$ are sampled based on equations \eqref{randomsample2}, respectively (lines 4-5). Based on our current time $t$ and the time to the next period, we can calculate the beginning of the next regulation period $t^{Start}$ (line 6). The deviations up until $t^{Start}$ are set to zero (lines 8-10). Starting from period $t^{Start}$ for $\tau^{\text{D}(+/-)}_{\omega}$ periods up to $t^{End}$, the deviations are set based on the average function $f^{(+/-)}$ and a random error term $\epsilon$ (lines 11-13). Next the current time is updated to the  $t^{End}$ (line 14). In this way, we move through the time horizon until we reach the end $|\mathcal{T}|$. The process is repeated for each scenario and once for upward and once for downward regulation scenarios. 

\begin{figure}
 \begin{subfigure}[t]{.49\textwidth}
	\centering
	\includegraphics[width=0.9\textwidth]{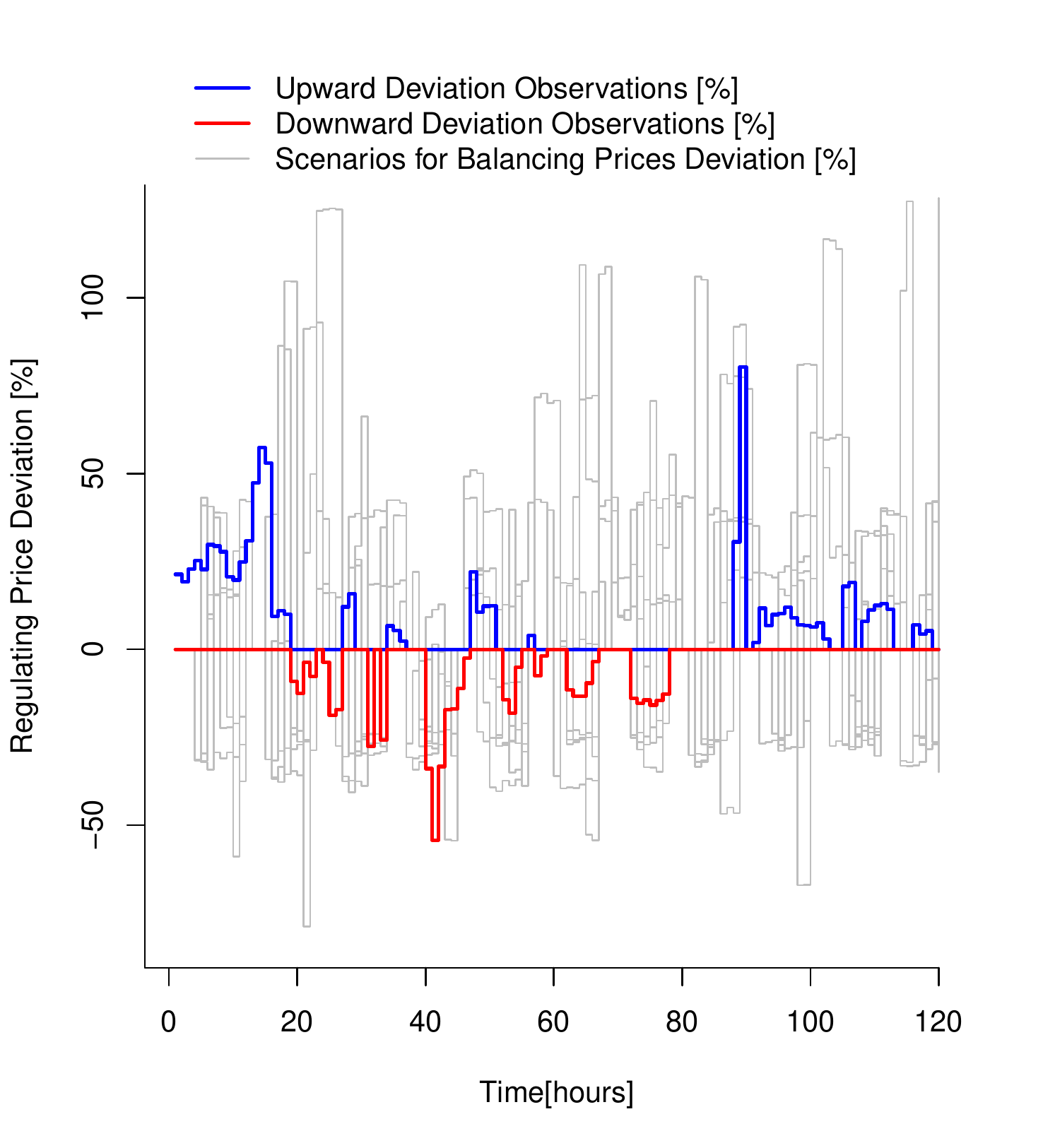}
	\caption{10 scenarios generated by Algorithm \ref{alg:balancingprices}.}
	\label{fig:BalancingScenariosOurs}
	\end{subfigure}
	 \begin{subfigure}[t]{.49\textwidth}
	\centering
	\includegraphics[width=0.9\textwidth]{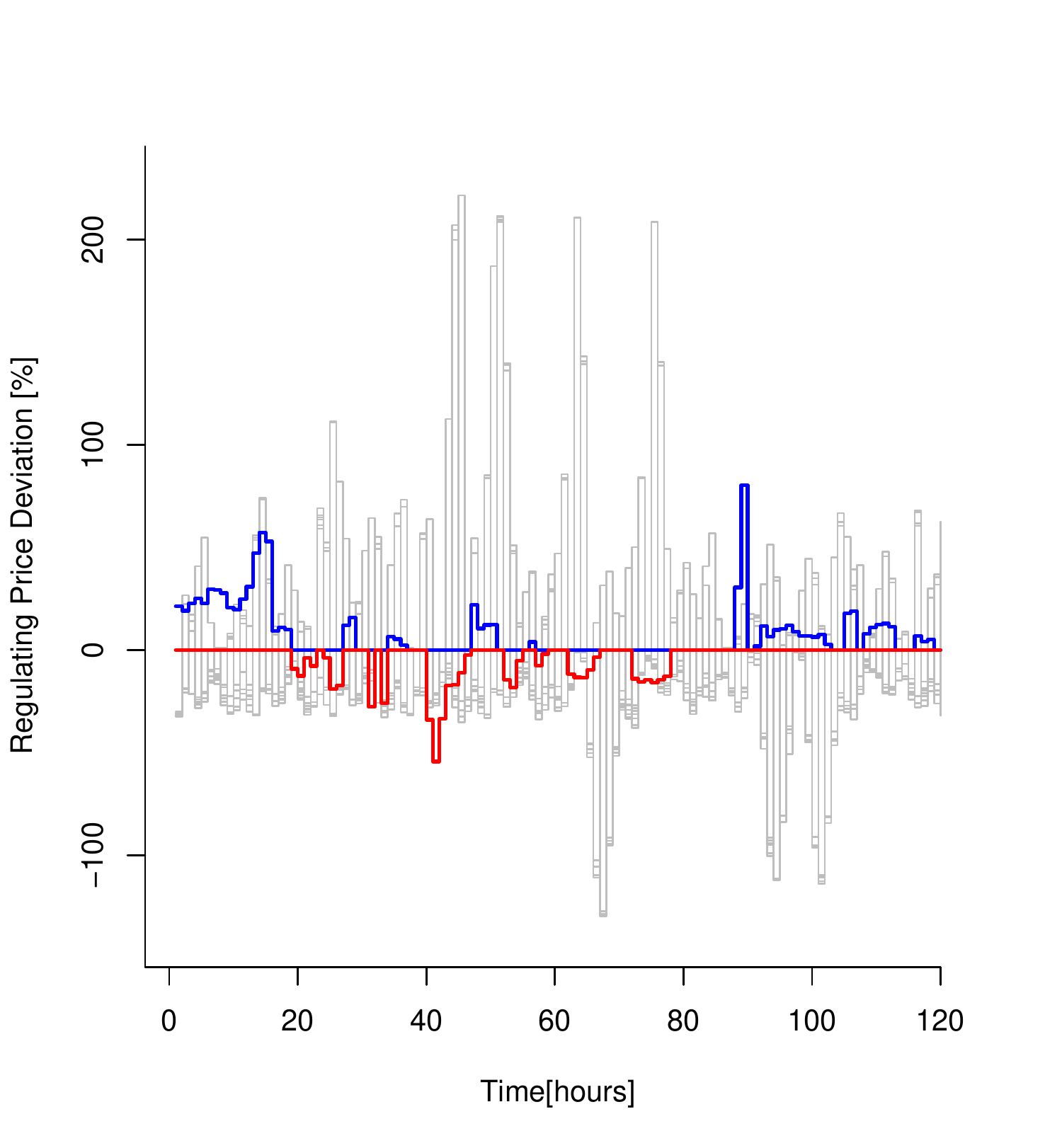}
	\caption{10 scenarios generated by the method proposed in \cite{olsson2008modeling}.}
	\label{fig:BalancingScenariosOlsson}
	\end{subfigure}
	\caption{Scenarios for balancing prices}
		\label{fig:BalancingScenarios}
\end{figure}

Since upward and downward regulation can not be activate at the same time, we calculate the final deviation scenario matrix as $\Delta \lambda_{t,\omega}=\Delta \lambda^{\text{UP}}_{t,\omega}-\Delta \lambda^{\text{DOWN}}_{t,\omega}$, where positive values of $\Delta \lambda_{t,\omega}$ represent upward regulation and the negative values downward regulation, respectively. 
Figure \ref{fig:BalancingScenariosOurs} shows a set of balancing prices scenarios generated by Algorithm \ref{alg:balancingprices} compared to the real observations. In comparison to scenarios generated by the method in \ref{fig:BalancingScenariosOlsson}, we can see the increased variability of regulating prices in the scenarios generated by Algorithm \ref{alg:balancingprices}. This is due to the fact that the prices are not based on time-series forecasts like in  \ref{fig:BalancingScenariosOlsson} but on the observed duration for upward and downward regulation periods.  To obtain the final prices the deviation value $\Delta \lambda_{t,\omega}$ is multiplied with respective day-ahead market price.


\section{Operational scheduling and bidding method}\label{sec:bidding}

The overall method, which allows the DH operator to schedule the production and determine the bidding curves for the day-ahead and balancing market, uses the two models presented in Section \ref{sec:model} with the scenarios generated by the methods in Section \ref{sec:uncertainty}. The optimization for one day in practice includes the following steps.

The day before the day in question, the \textit{day-ahead market optimization} \eqref{eq:dayahead1}-\eqref{eq:dayahead3} is solved as two-stage stochastic programming. The model includes scenarios representing the uncertainty regarding day-ahead market electricity prices (Section \ref{sec:scenarios_dayahead}), wind power production (Section \ref{sec:scenarios_wind}) and solar heat production (Section \ref{sec:scenarios_solar}) for at least 24 hours. The scenarios are generated using the Monte Carlo simulation and clustering technique described in Section \ref{sec:scenario_generation}. The planning horizon can be considered as longer than 24 hours in a rolling horizon manner to include future days into the optimization to get better approximation of the thermal storage behaviour, which can store heat longer than just 24 hours. The optimal values of the variables $p^{\text{BID}}_{t,\omega}$ in \eqref{eq:dayahead1}-\eqref{eq:dayahead3} return the bidding amounts for each hour $t \in \lbrace 1,\ldots, 24\rbrace$, while each scenario $\omega$ sets one step in the bidding curve. As  constraints \eqref{equalitiybids} and \eqref{nondecreasingbids} ensure the same production amounts for the same electricity prices and increasing production amounts for increasing prices, the optimal values $p^{\text{BID}}_{t,\omega}$ result automatically in a non-decreasing step-wise bidding curve. The bidding prices for each step in the bidding curve are the respective electricity price forecast values $\lambda_{t,\omega}$.

After the day-ahead market is cleared, the real electricity prices for each hour become available and the won bids can be determined (i.e. the hours where the bidding price was equal or below the market price). In hours with won bids, the DH operator is committed to provide the offered amount of power, otherwise the caused imbalance is penalized with a payment. However, imbalances from other operators on the market offer an opportunity for profit. The balancing market is used by the TSO to reduce the imbalances in the system by accepting new bids for additional power or reducing production. Thus, we can use the flexibility in our portfolio of production units to also offer upward and downward regulations bids in the balancing market. As the \textit{balancing market} has a time horizon of only one hour and is closed shortly before this hour, an optimization needs to take place every hour before the balancing market closes. Model \eqref{eq:balancing1} to \eqref{eq:balancing3} optimizes the production for the next hour taking the committed production from the day-ahead market into account. Furthermore, the model can take several hours into the future into account to anticipate impact on the remaining hours of the day. The model is again a two-stage stochastic program considering the balancing market price scenarios (see Section \ref{sec:scenarios_balancing}) for all hours and wind power scenarios for later on the day (we assume that the wind power for the next hour can be predicted accurately). Again, the optimal values of $p^{\text{UP}}_{t,\omega}$ and $p^{\text{DOWN}}_{t,\omega}$ result automatically in a non-decreasing or non-increasing step-wise bidding curves representing upward and downward regulation bids, respectively. The bidding prices for each step in the bidding curve are the respective electricity price forecast values $\lambda^\text{UP}_{t,\omega}$ and $\lambda^\text{DOWN}_{t,\omega}$. This step is repeated by the operator for each hour.

\section{Case study}\label{sec:casestudies} 

We use the Hvide Sande district heating system\footnote{see Hvide Sande Fjernvarme A.m.b.A., \href{https://www.hsfv.dk/}{https://www.hsfv.dk/}} in Western Jutland, Denmark, as a case study to evaluate our method. However, the method presented in this paper is applicable to all district heating systems with a portfolio of units, because the models in Section \ref{sec:model} are formulated in a general manner and the scenario generation methods \ref{sec:uncertainty} can be replaced by other available forecasting techniques without changing the overall methodology. 

\begin{figure}
	\centering
	\includegraphics[width=0.55\textwidth]{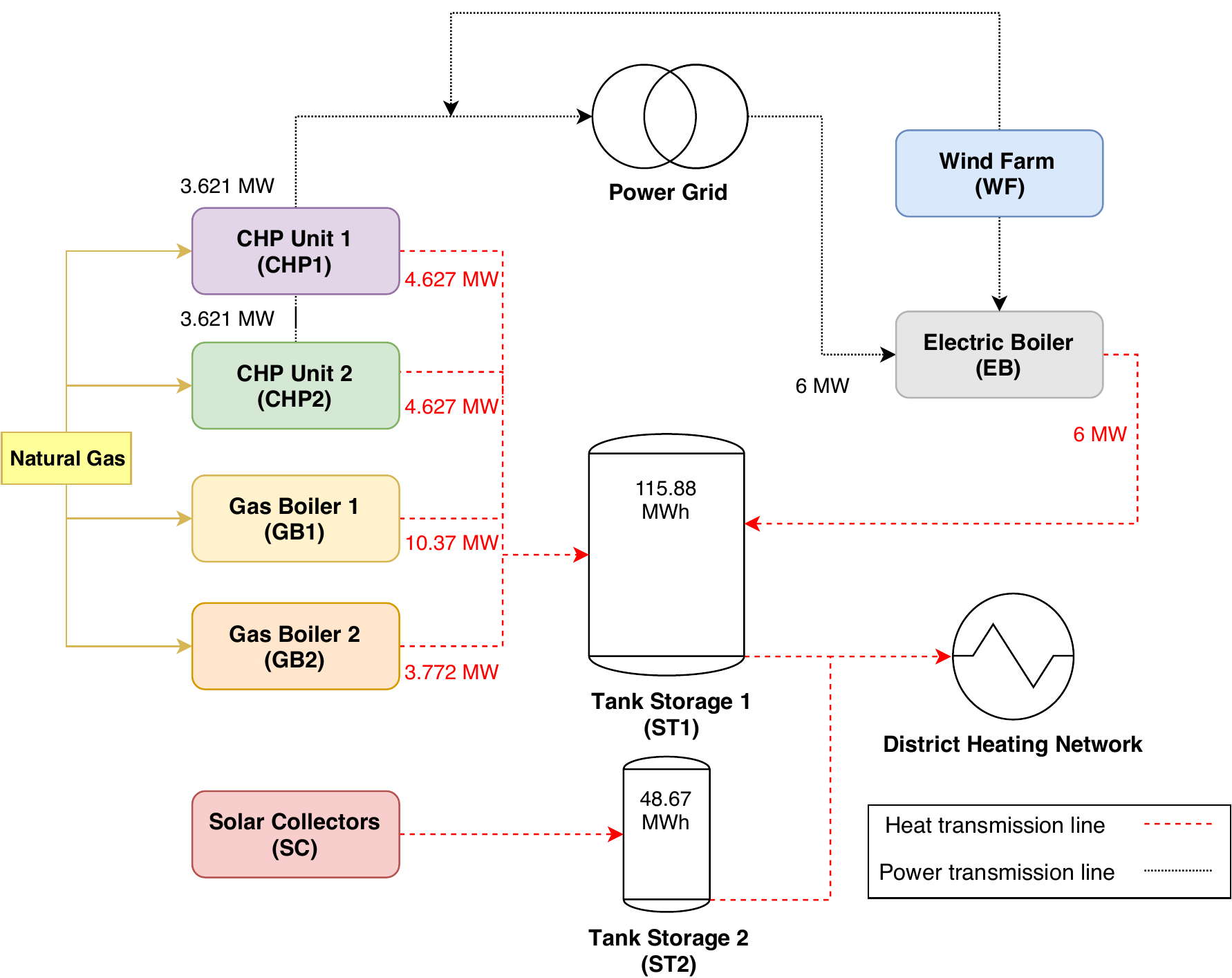}
	\caption{Flowchart of the Hvide Sande district heating system}
	\label{HvideSande_Flowchart}
\end{figure}

An overview of the Hvide Sande system is given in Figure \ref{HvideSande_Flowchart}. It has two small gas-fired CHP units (CHP1 and CHP2) acting on the electricity market and feeding heat to the district heating system as well as two gas boilers (GB1 and GB2) units with dispatchable heat production. Stochastic renewable heat production comes from a solar collector field (SC), which is considered as one unit. Finally, it is also possible to produce heat from electricity using an electric boiler (EB). The electricity can be bought from the power grid as a regular consumer or using a special tariff. This tariff consists of a tax benefit for operating the electric boiler, in which the amount of power injected by the own wind farm (WF) into the grid is at the same time consumed by the electric boiler. This synchronous operation of both units help the power system to reduce imbalances and provides cheap heat production.  
The DH system has two thermal storages, where one (ST1) is connected only to the solar collector field and the second storage (ST2) is used by all other units. The parameters for costs and capacities as well as the connections between units are given in Table \ref{TechnicalDataandcosts}. Furthermore, the table shows to which set the units belong.

\begin{table}
	\centering
	\footnotesize
	\caption{Characteristics of the production units and thermal storages}
\begin{tabular}{llrrrrrrrr}\toprule

      Unit &        Set & $C_{u}^{H}$ & $C_{u}^{T}$ & $\overline{Q}_{u}$ & $\overline{P}_{u}$ & $\varphi_{u}$ & $A^{\text{DH}}_{u}$ & \multicolumn{2}{c}{$A^{\text{S}}_{u,s}$}             \\ \cmidrule(l){9-10}

           &            &            &            &            &            &            &            &         ST1 &         ST2 \\\midrule

      CHP1 & $\mathcal{U}^{\text{CHP}}$ &     689.01 &          - &       4.63 &       3.62 &       1.28 &          0 &          1 &          0 \\

      CHP2 & $\mathcal{U}^{\text{CHP}}$ &     689.01 &          - &       4.63 &       3.62 &       1.28 &          0 &          1 &          0 \\

       GB1 & $\mathcal{U}^{\text{H}}$ &     401.30 &          - &      10.37 &       0.00 &          - &          0 &          1 &          0 \\

       GB2 & $\mathcal{U}^{\text{H}}$ &     416.29 &          - &       3.77 &       0.00 &          - &          0 &          1 &          0 \\

        EB & $\mathcal{U}^{EL}$ &     359.98 &      49.52 &       6.00 &       0.00 &       1.00 &          0 &          1 &          0 \\
  SC & $\mathcal{U}^{\text{RES}}$ &       0.00 &          - &     100.00 &       0.00 &          - &          0 &          0 &          1 \\
        WF & $\mathcal{G}$ &       0.00 &          - &       0.00 &       0.00 &          - &          - &          - &          - \\

    \cmidrule(l){3-10} \cmidrule(l){3-10}

           &            & $\overline{S}$ & $\underline{S}$ & $\sigma_0$ &            &            &            &            &            \\ \cmidrule(l){3-10}

       ST1 & $\mathcal{S}$ &     115.88 &      0.00 &      57.94 &            &            &            &            &            \\

       ST2 & $\mathcal{S}$ &      48.67 &       0.00 &      24.34 &            &            &            &            &            \\
        \bottomrule
\end{tabular}

	\label{TechnicalDataandcosts}
\end{table}

\section{Analysis of experimental results}\label{sec:results}

To evaluate our approach, we have to determine the real costs and behaviour of the system. The actual wind power production, solar thermal production and heat demand values are obtained from the Hvide Sande district heating system for the year 2017. The day-ahead, upward and downward electricity prices are taken from the  NordPool market for the bidding area DK1 (where Hvide Sande is located). This data is public and can be downloaded from \cite{WelcomeW75:online}. The data basis for forecasting and scenario generation is historical data from 15 days before the day in question. The input data for wind speed, solar radiation and ambient temperature are randomly perturbed values of the real data. The overall evaluation process includes the following steps:

\begin{enumerate}
    \item Before day-ahead market closure for day $d$ (Day $d-1$): Create scenarios for the day-ahead market optimization and solve optimization model \eqref{eq:dayahead1}-\eqref{eq:dayahead3} using thermal storage level from the day before. Submit bids to the day-ahead market.
    \item After day-ahead market closure for day $d$ (Day $d-1$): Evaluate the day-market bids with the now known electricity prices and save production amounts of won bids.
    \item Each hour on day $d$:
    \begin{enumerate}
        \item Before the closure of the balancing market at hour $t$ on day $d$: Create scenarios for the balancing market optimization, include the committed power production amounts from the day-ahead market and solve optimization model \eqref{eq:balancing1}-\eqref{eq:balancing3}. 
        \item Evaluate the balancing-market bids with the now known balancing electricity prices, fix the committed production amounts and resolve the model to get actual costs and thermal storage levels.
    \end{enumerate}
    \item Move to the next day
\end{enumerate}

The forecasting and scenario generation are implemented in R 3.2.2, while the optimization models are built
in GAMS 24.9.2 using CPLEX 12.1.1 to solve them. All experiments were executed on the DTU HPC Cluster using 2xIntel Xeon Processor X5550 and 24 GB memory RAM. For the results presented in the remainder of this section, we use a rolling horizon of three days in the day-ahead optimization problem and 12 hours for the balancing market problem. To correlate different scenarios of RES with electricity prices, we generate $n$ different scenarios for wind power and solar heat production and $m$ scenarios for electricity prices. The combination of all scenarios results in a total number of $m \times n$ scenarios. For the sake of simplicity we generate $n=10$ scenarios of RES production for the experiments that consider bidding curves. 

\subsection{Influence of uncertainty and number of bidding curve steps on the day-ahead market results}

\begin{figure}[tbp]
\centering
\includegraphics[width=0.6\textwidth]{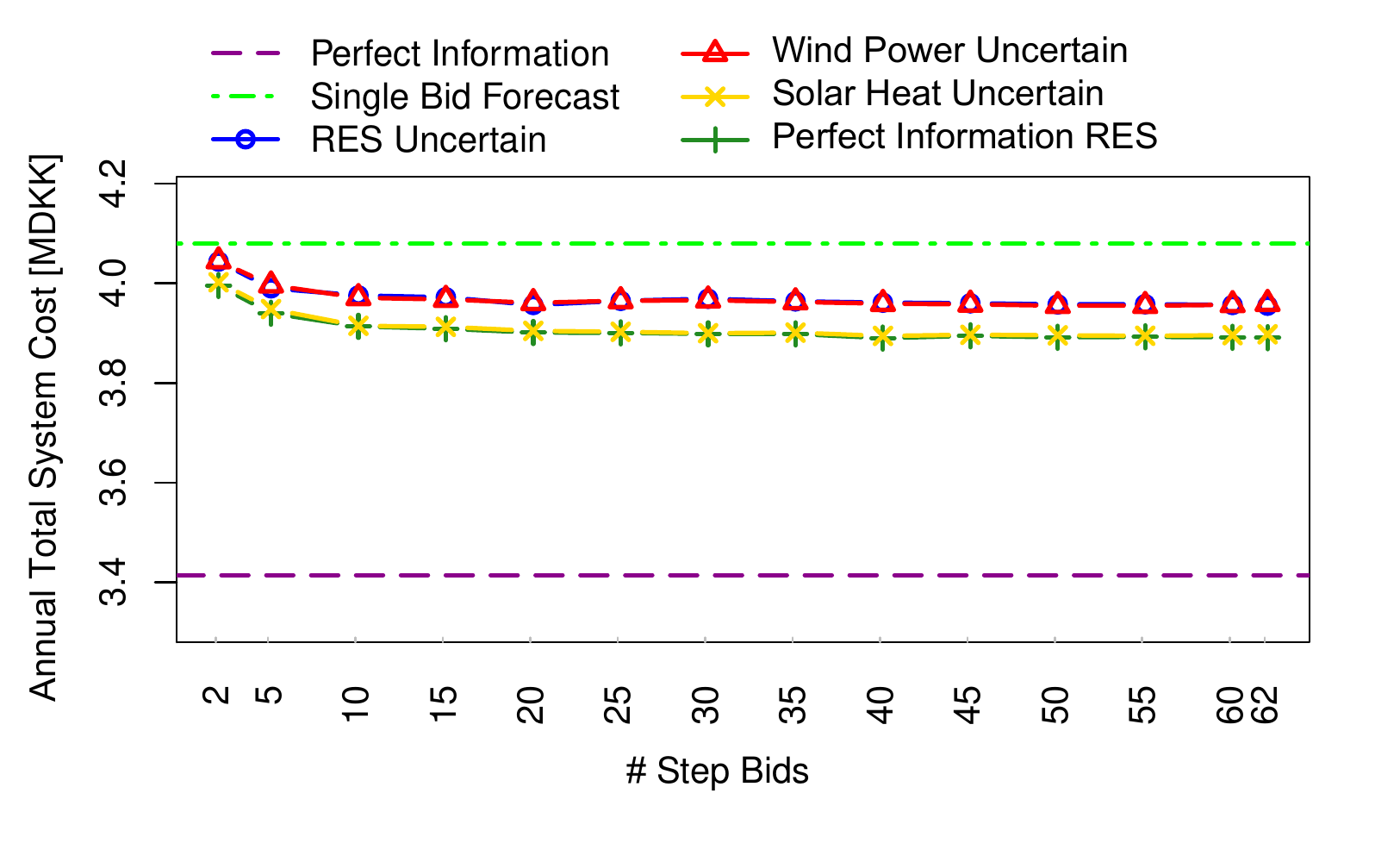}
\caption{Comparison of different uncertainty setups and number of steps in the bidding curves in the day-ahead market optimization. The values shown are total annual system cost. }
\label{fig:results_uncertainty_steps}
\end{figure}

In the first experiment, we concentrate on the bidding results from the day-ahead market optimization problem \eqref{eq:dayahead1}-\eqref{eq:dayahead3} only. We compare the total annual costs of different setups regarding uncertainty consideration, i.e., which values are known or unknown, and the number of electricity price scenarios resulting in the steps for the bidding curves. The results are given in Figure \ref{fig:results_uncertainty_steps}, where the x-axis represents the number of steps in the bidding curve (i.e. the number of electricity price scenarios) and the y-axis represents the total annual system costs. The depicted lines show the results of different setups regarding uncertainty consideration. The theoretically best result is given by considering that we have perfect knowledge about the future electricity prices and RES production (\textit{Perfect Information}). However, this value can never be reached in practice due to the uncertainty and, therefore, serves only as benchmark. Another value to compare to is a bidding method that submits bids according to the expected electricity price (\textit{Singe Bid Forecast}), i.e., the model considers no electricity price scenarios but the expected value resulting in one bidding amount and price for each hour. This approach is often used in practice. The other four approaches consider the model from Section \ref{sec:dayaheadmodel} to create bidding curves based on uncertain electricity prices. We compare four cases regarding the information about RES production: scenarios for wind power and solar thermal production (\textit{RES Uncertain}), scenarios for wind power and perfect information about solar thermal production (\textit{Wind Power Uncertain}), scenarios for solar thermal and perfect information about wind power production (\textit{Solar Heat Uncertain}) as well as perfect information regarding both (\textit{Perfect Information RES}). 

The results from Figure \ref{fig:results_uncertainty_steps} indicate that considering the solar thermal production as uncertain and modelling it as scenarios does not deteriorate the costs significantly compared to the case where the RES production is known. On the other hand, considering the wind power production as uncertainty captured in  scenarios, has an impact on the costs and leads to an increase in the cost of approx. 62000 DKK. Similar results are achieved when considering both RES production sources as uncertain. Based on this results, we can conclude that especially the uncertainty of the wind power production has an influence on the systems costs. This behaviour can be explained based on the fact that the wind power production has a direct effect on the power amount that is traded on the electricity market and therefore on the profits obtained. In contrast, the uncertainty of solar thermal production has no large effect due to the thermal storage in the system, which smooths the effect on the heat production and therefore also the costs. The factor that has the greatest impact on the operational cost is not having  information about the day-ahead prices (Perfect Information). In this case having perfect information of RES and uncertain day-ahead electricity prices increased the annual system cost by approx. 500,000 DKK (around 12.5\% of the total system cost). However, under the real-world condition that RES and electricity prices are uncertain, using stochastic programming to generate bidding curves decreases the cost by ca. 120,000 DKK per year (3\% of the total system cost) compared to the \textit{Single Bid Forecast}. 

 Figure \ref{fig:results_uncertainty_steps} also shows the influence of the number of steps in the bidding curves. For this experiment the number of clusters in the PAM algorithm was varied (see Section \ref{sec:scenario_generation}) to obtain different numbers of scenarios representing the number of steps in the bidding curve. We compare in total 14 scenario set sizes ranging from 2 to 62 scenarios, which are the minimum and the maximum number of steps allowed to submit to the NordPool market \cite{nordpool43:online}, respectively. The results show a reduction of costs when the number of steps is increased from two to 20 steps. In this case, including more steps does not lead to further significant reductions in costs. 

Based on the analysis in this section, we can conclude that using bidding curves, in particular with at least 20 steps, created from our stochastic program can reduce the annual system cost in particular compared to single bids based on price forecasts. Furthermore, the uncertainty of wind power production influences the results more than the uncertainty regarding solar thermal production.

\subsection{Impact of special tariff for the electric boiler}

\begin{figure}
\centering
\includegraphics[width=0.6\textwidth]{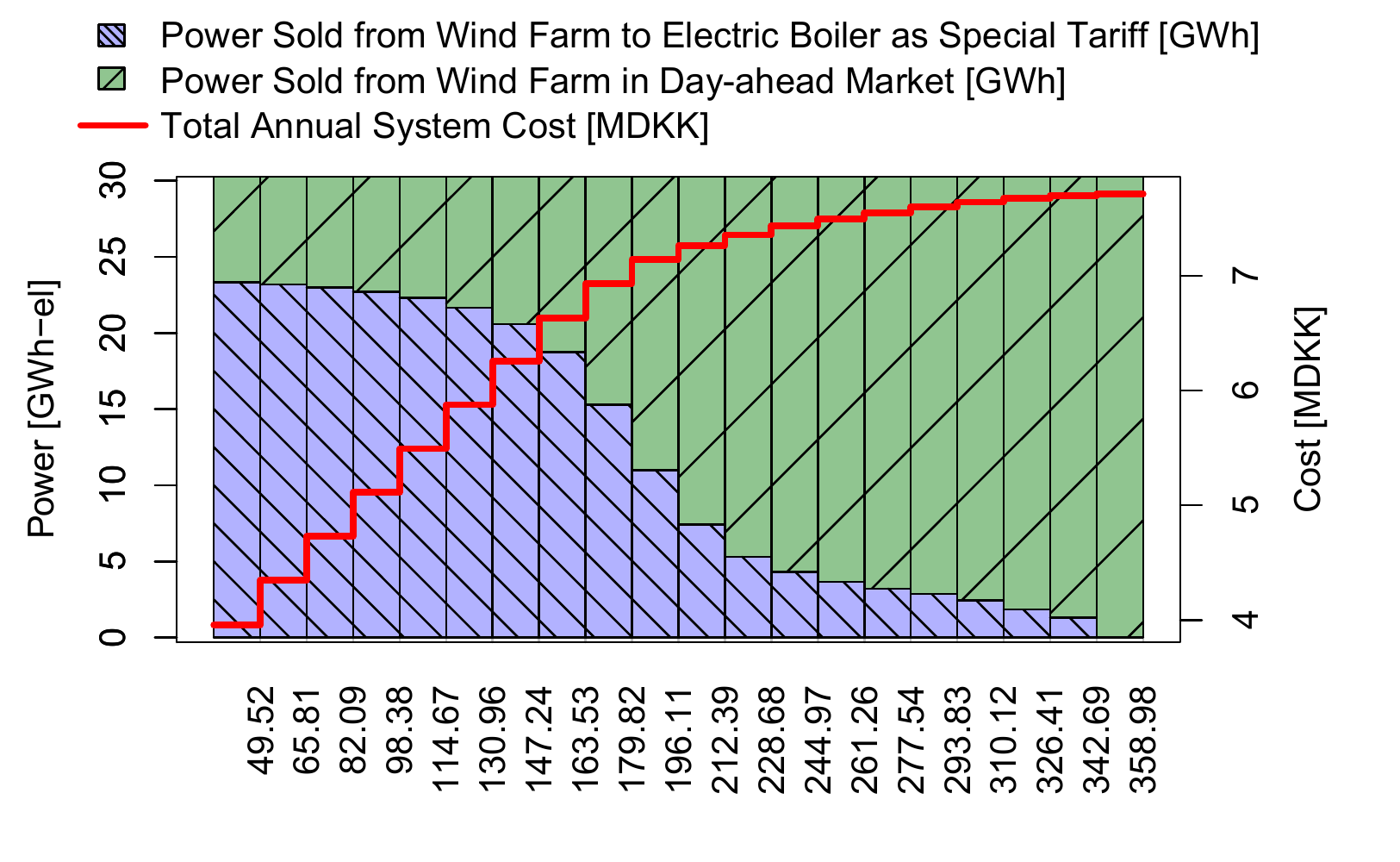}
\caption{Power from the wind turbines used for the electric boiler as special tariff or traded on the day-ahead market as well as total annual system cost. The values are given for varying special tariff operation cost of the electric boiler.}
\label{fig:results_elboilertariff}
\end{figure}

As mentioned in the problem description in Section \ref{sec:model}, we assume a special tariff (in terms of tax reduction) if the electric boiler is "using" power that we provide with our wind farm. In this section, we want to analyze the influence of this tariff on the trading on the day-ahead market. The operational cost under the special tariff were given with 49.52 DKK/MWh-heat. Figure \ref{fig:results_elboilertariff} shows the impact on the annual system cost and share of wind power used for the electric boiler and traded on the day-ahead market, respectively, when the tariff is increased in equal step sizes up to the normal operational cost (when fed with power from the grid without special tariff). 

Figure \ref{fig:results_elboilertariff} shows clearly the benefits from having a special agreement when feeding in wind power and therefore receiving a special tariff on the electricity consumption. First, the total annual system cost drastically increase when the special tariff gets more expensive. This is obvious as the production of heat from electricity is getting more expensive. Furthermore, it can be seen that the amount of wind power traded on the day-ahead market increases with a higher tariff, because the income from the market is more promising in most of the hours in the year. This means, using the special tariff for the electric boiler is only beneficial, if the income from the market is expected to be less than the benefit from using the wind power for the electric boiler. This margin is getting smaller with increasing special tariff, resulting in higher trade volumes on the day-ahead market.

This results indicate that DH operators can greatly benefit from receiving a special agreement with respect to using own RES power generation for heat production.

\subsection{Analysis of yearly production}

\begin{figure}
\begin{subfigure}{1\textwidth}
\centering
\includegraphics[width=0.8\textwidth]{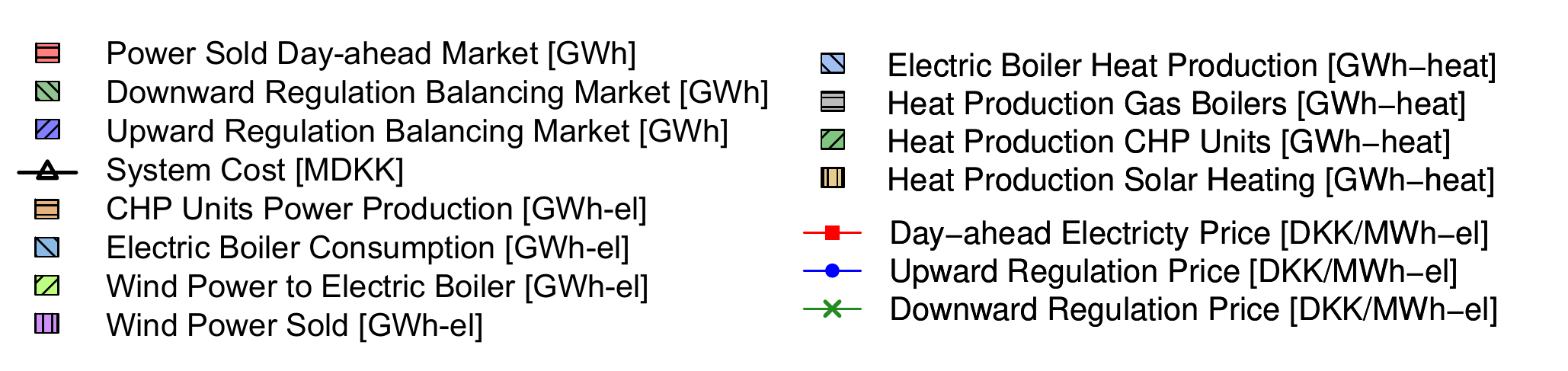}
\caption{Legend}
\label{fig:legend}
\end{subfigure}
	\begin{subfigure}[t]{.50\textwidth}
	\centering
	\includegraphics[width=\textwidth]{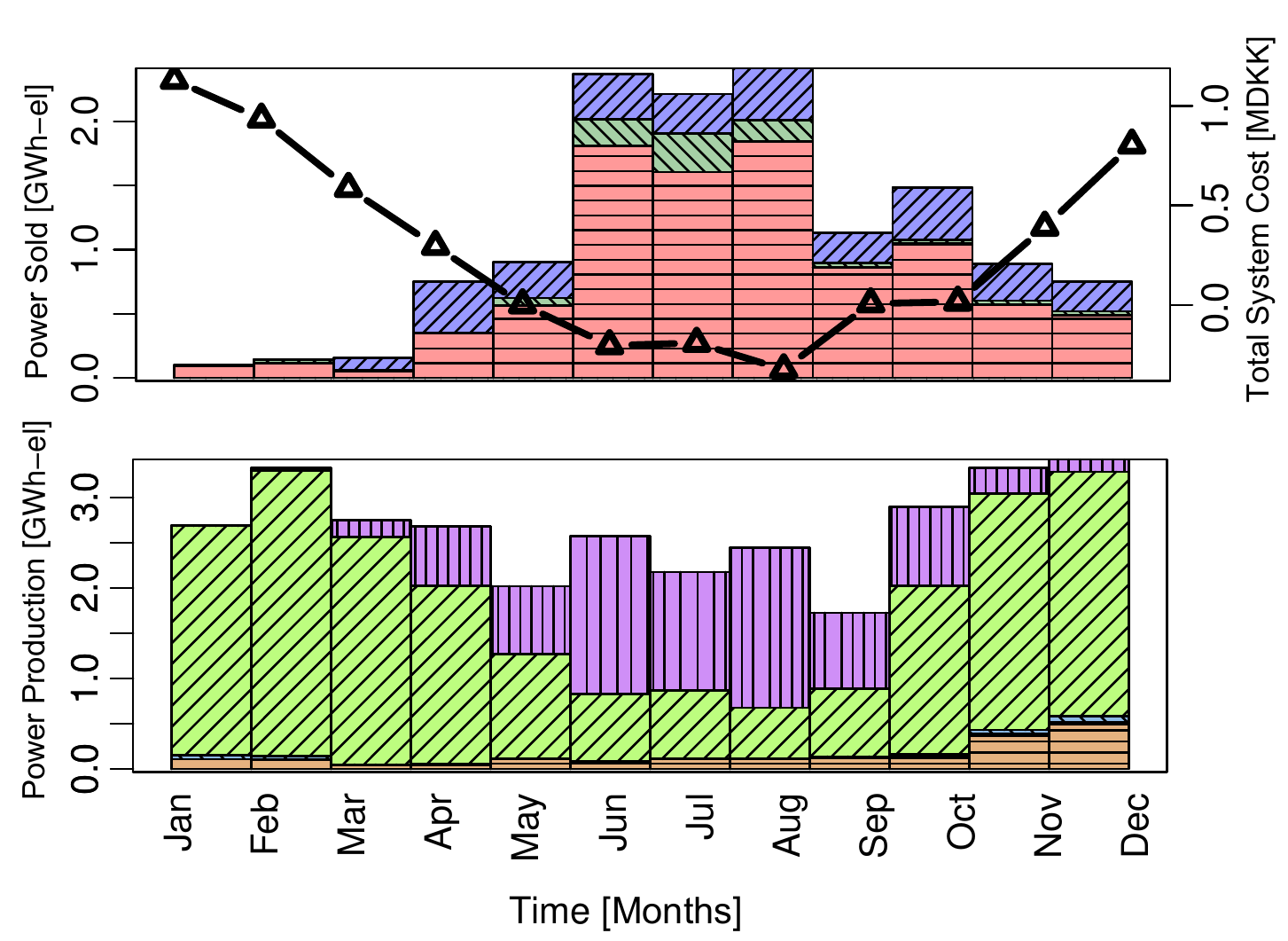}
	\caption{Power sold on the markets and system cost (top) as well as power production by the different units (bottom) }
	\label{fig:results_power_annual}
	\end{subfigure}
	\quad
   \begin{subfigure}[t]{.48\textwidth}
	\centering
	\includegraphics[width=\textwidth]{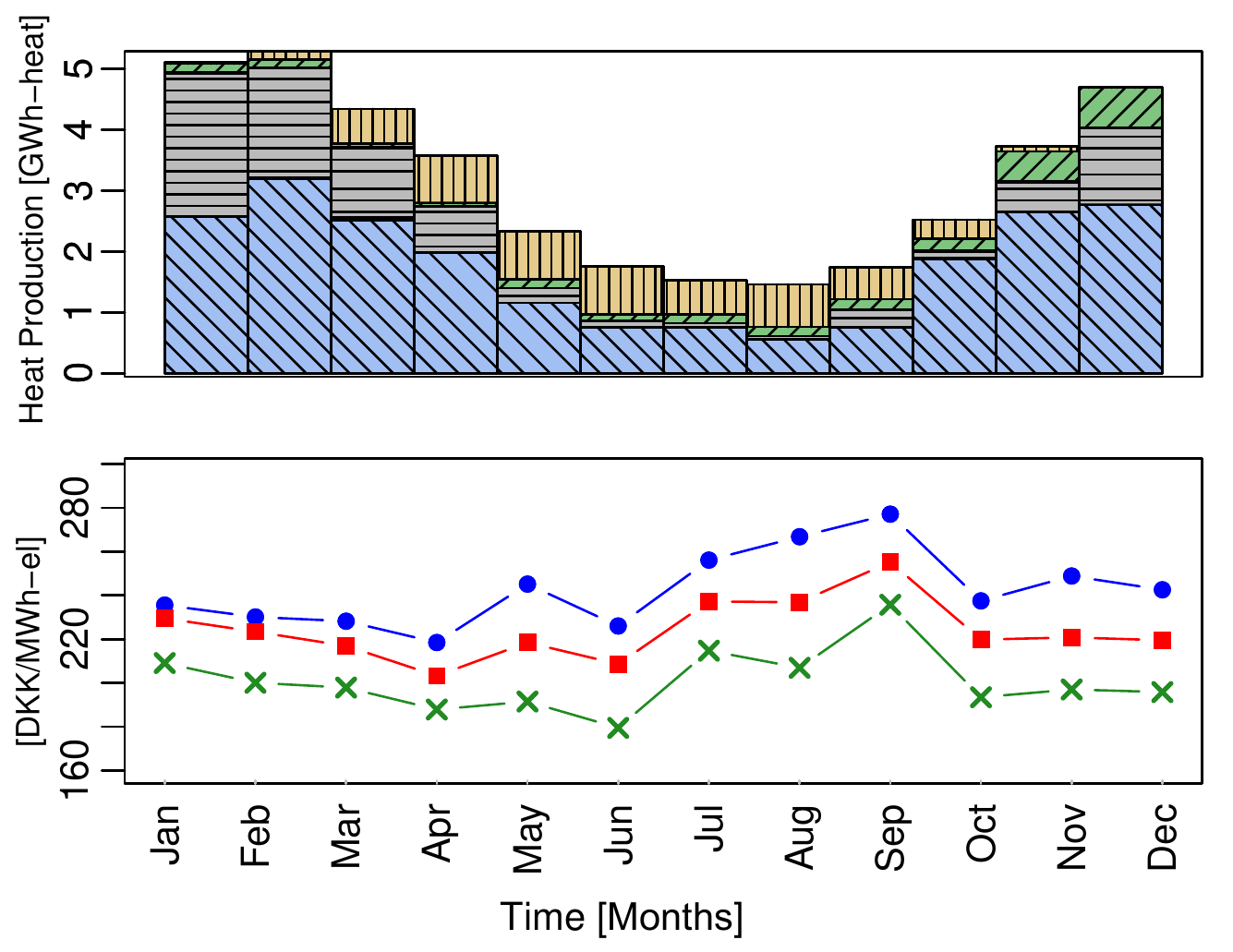}
	\caption{Heat production by the different units (top) and average electricity prices (bottom)}
	\label{fig:results_heat_annual}
	\end{subfigure}
	\caption{Annual system behaviour on a monthly time-scale} \label{fig:results_annual}
\end{figure}

In this section we provide the results of the annual system behaviour when using the solution approach for day-ahead market and balancing market optimization presented in Section \ref{sec:bidding}. The results and values for power production and trading, heat production and electricity prices are consolidated on a monthly basis in Figure \ref{fig:results_annual}. The legend can be found in Figure \ref{fig:legend}.

Figure \ref{fig:results_power_annual} (top) shows the monthly system cost and the amounts of power traded on the day-ahead market as well as down regulation bought and upward regulation sold on the balancing market. One observation from this figure is that the monthly costs are significantly lower during the summer period due to two reasons. First, the amount of power traded on the different electricity markets is higher during the summer resulting in higher profits. Also the electricity prices are slightly higher during the summer (see Figure \ref{fig:results_heat_annual} (bottom)). Second, the heat demand is lower during the summer resulting in lower total costs (see Figure \ref{fig:results_heat_annual} (top)). Furthermore, from Figure \ref{fig:results_heat_annual} (top) it can be seen that the solar thermal production is higher during summer resulting in less heat needed from the more expensive other units. 

A second observation is that the trades on the balancing markets have a higher volume during summer and fall. This behaviour can be explained by taking the power production on a unit-basis into account as provided in Figure \ref{fig:results_power_annual} (bottom). During the summer month less of the power is used for heat production, because there is a lower heat demand, and therefore available for trading on the electricity markets. 

The results show that the optimization method makes use of the fact that the units are considered as one portfolio and thereby having a flexibility with respect to the production. The trading and production behaviour adjusts itself to the best combination in the different seasons to get lowest heat production costs and highest incomes from the markets. The specific daily system behaviour in case of regulation activities is further analyzed in Section \ref{sec:results_upwarddownward}.

\subsection{Value of including balancing market trading}

\begin{table}
\footnotesize
\centering
\begin{tabular}{lrr}
\toprule
  Setting  &  Total System Cost [DKK] &    $\Delta$        \\\midrule

Perfect information incl. balancing market  &  2,499,205 &    -        \\

Perfect information excl. balancing market  &  3,414,310 &   +37\%          \\

Stochastic incl. balancing market  &  3,655,798 &   +7\%          \\

Stochastic excl. balancing market  &  3,956,530 &    +8\%        \\
        \bottomrule
\end{tabular}
\caption{Comparison of annual system cost in different setups of the solution approach}
\label{tab:results_includebalancing}
\end{table}

The next analysis investigates the value of including the balancing market trading into the solution approach. Therefore, we compare two settings: Using the solution approach from Section \ref{sec:bidding} with and without the balancing market optimization. Furthermore, we run both settings once with perfect information about electricity prices, wind and solar production and once in a stochastic programming setting with scenarios (as presented in the model formulation in Section \ref{sec:model}). 

The total annual system cost for those four cases in Table \ref{tab:results_includebalancing} show that even if we have perfect information about the future ignoring trading on the balancing market will increase the total system cost immensely, in this case by 37\%. This indicates that a high degree of income can be obtained from the balancing market. The results with perfect information are theoretical values as those can not be reached in a real world application due to the uncertainty at the time of planning. This means that when modelling the uncertainty regarding prices and production in a stochastic setting, the system cost naturally increase, here by 7\%. However, lower cost are still achieved by including the balancing market as a second step to avoid imbalances and another opportunity for trading. Not considering the balancing market leads here to 8\% higher system cost for the entire year.

\begin{figure}
\centering
\includegraphics[width=0.65\textwidth]{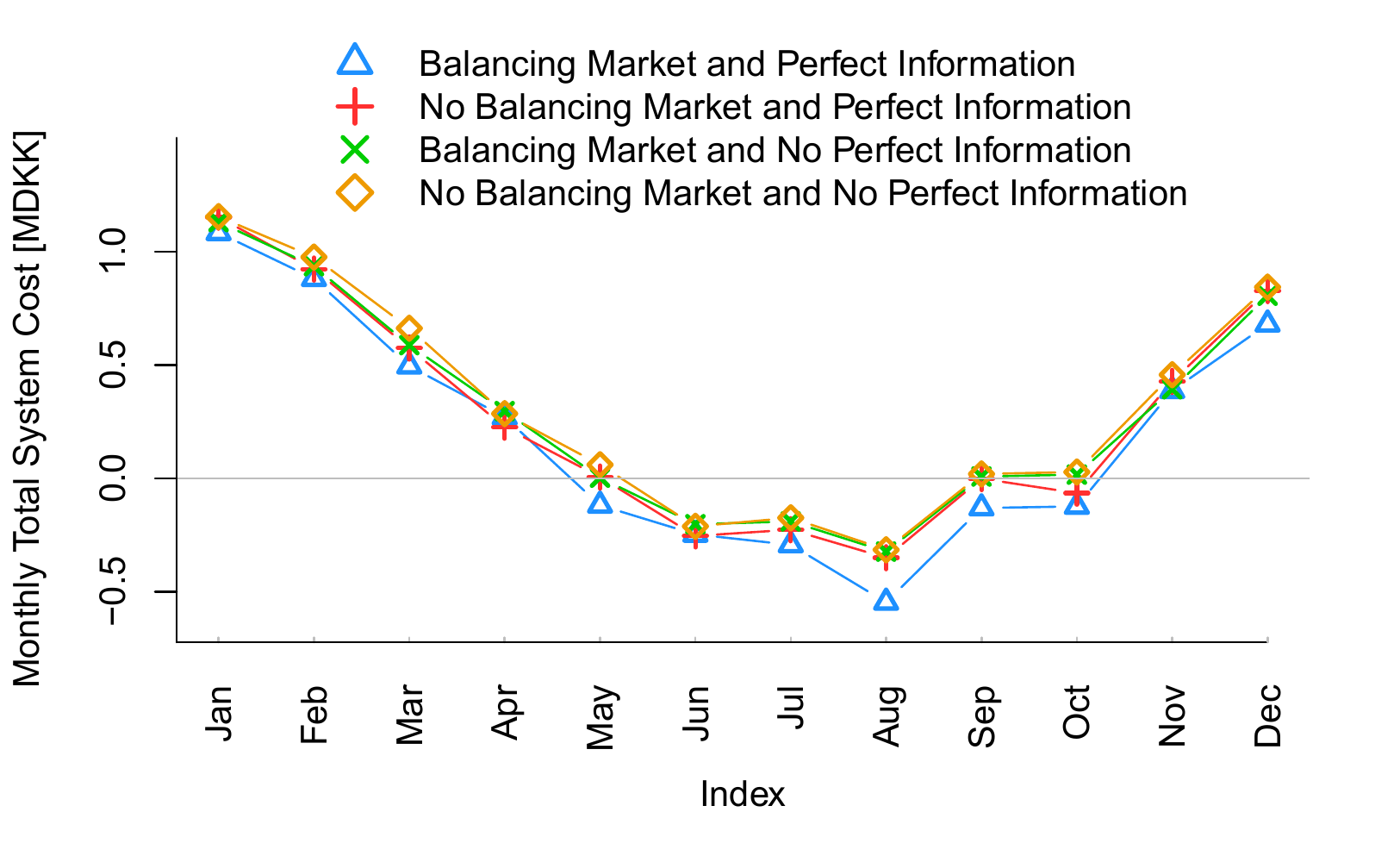}
\caption{Comparison of monthly system cost in different setups of the solution approach}
\label{fig:results_includebalancingmontly}
\end{figure}

These conclusions are mostly independent from the actual months or seasons as it can be seen from Figure \ref{fig:results_includebalancingmontly}. Here the monthly cost for the four settings follow a similar ranking in each month. 

\subsection{Behaviour of system in case of upward and downward regulation}\label{sec:results_upwarddownward}

To further investigate the benefits of trading in the balancing market, we analyze the obtained production schedules for four representative days of the year where upward and downward regulation was offered. Section \ref{sec:upwardprices} and \ref{sec:downwardprices} each analyze two specific days in which upward and downward regulation was provided, respectively. The legend used for the production schedule figures in this section is the same as in Figure  \ref{fig:legend}.

\subsubsection{Upward regulation}  \label{sec:upwardprices}

\begin{figure}
	\begin{subfigure}[t]{.44\textwidth}
	\centering
	\includegraphics[width=1.0\textwidth]{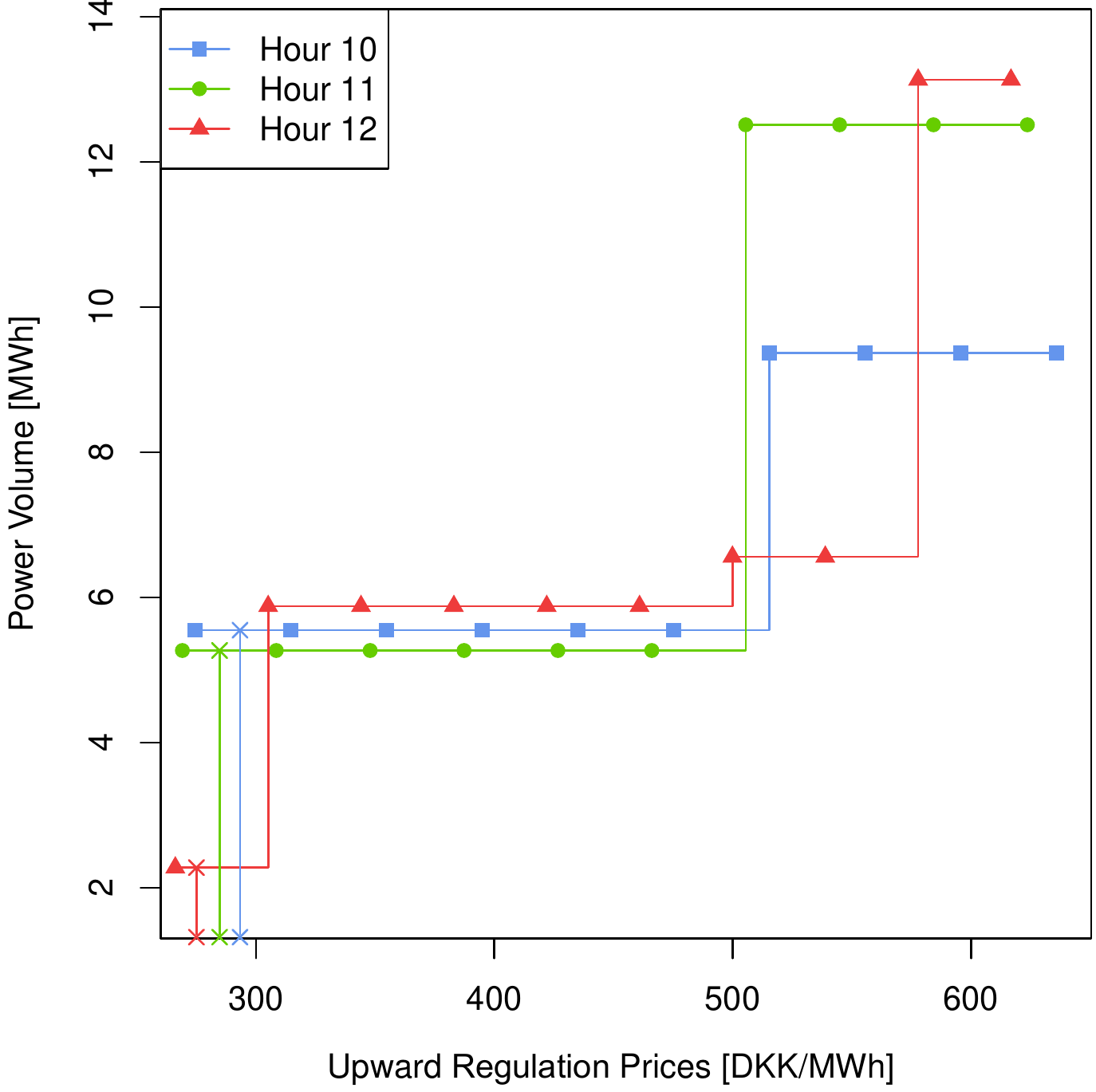}
	\caption{Bidding curves and won bids (marked with $\times$)}
	\label{fig:Upcase1BiddingCurves}
	\end{subfigure}
   \begin{subfigure}[t]{.54\textwidth}
	\centering
	\includegraphics[width=1.0\textwidth]{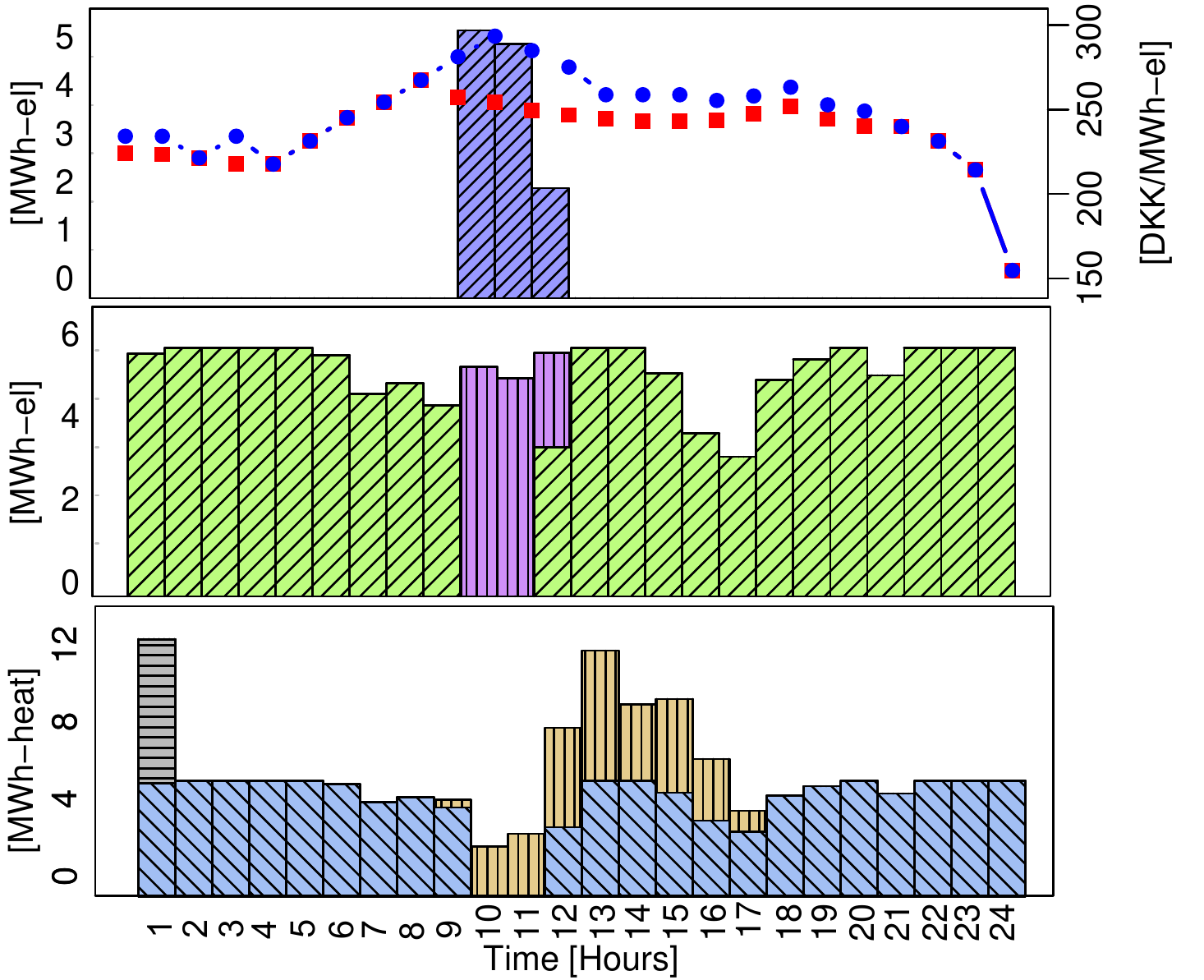}
	\caption{System operation: upward regulation amount including day-ahead and balancing prices (top), power production (middle), heat production (bottom)}
	\label{fig:Upcase1Operation}
	\end{subfigure}
	\caption{Upward regulation provided on 9th March 2017 } \label{fig:UpCase1}
	\end{figure}

The first case for upward regulation is presented in Figure \ref{fig:UpCase1}. Figure \ref{fig:Upcase1BiddingCurves} shows the bidding curves for the hours in which upward regulation was won by the DH operator. The vertical lines delimited with "$\times$" represent the real upward prices for those hours and the corresponding power production offered at such prices. Figure \ref{fig:Upcase1Operation} shows the system behaviour and is divided into three parts: 
upward regulation volume per hour including prices (top), hourly power production per unit (middle) and hourly heat production per unit (bottom). As we can see from  \ref{fig:Upcase1Operation} (middle), the upward regulation in this case is entirely provided by the wind farm. Since no wind power was sold on the day-ahead market, the producer decides to bid the entire production of the wind farm into the balancing market for hours 10 and 11. In hour 12, the needed power volume for upward regulation is lower than the actual production from wind. Therefore, the remaining power is used to feed the electric boiler. This behaviour is confirmed by the heat production (Figure \ref{fig:Upcase1Operation} (bottom)). In hours 10 and 11 there is no production with the electric boiler but in hour 12.  

\begin{figure}
		\begin{subfigure}[t]{.45\textwidth}
	\centering
	\includegraphics[width=1.0\textwidth]{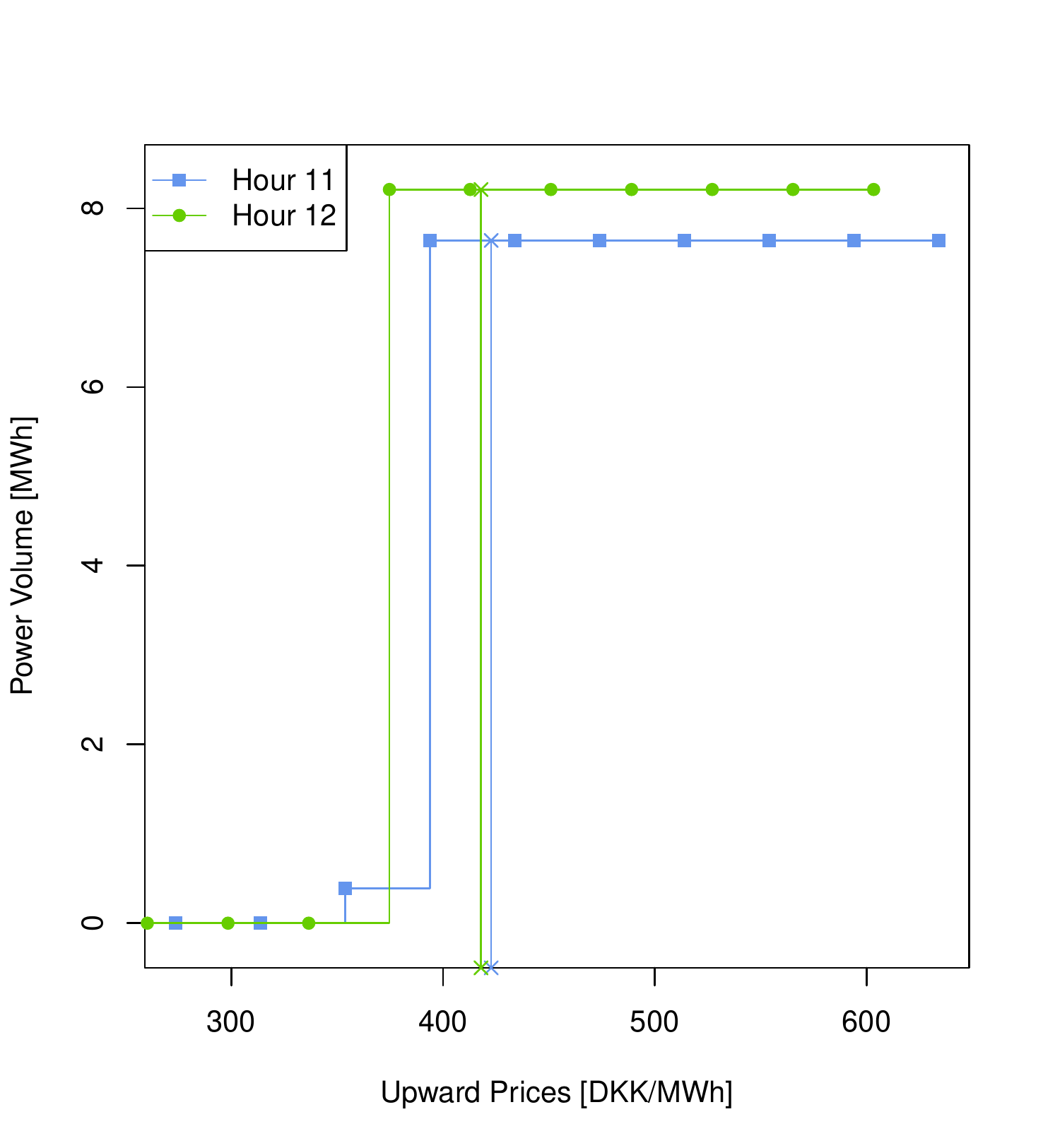}
	\caption{Bidding curves and won bids (marked with $\times$)}
	\label{fig:Upcase2BiddingCurves}
	\end{subfigure}
   \begin{subfigure}[t]{.52\textwidth}
	\centering
	\includegraphics[width=1.0\textwidth]{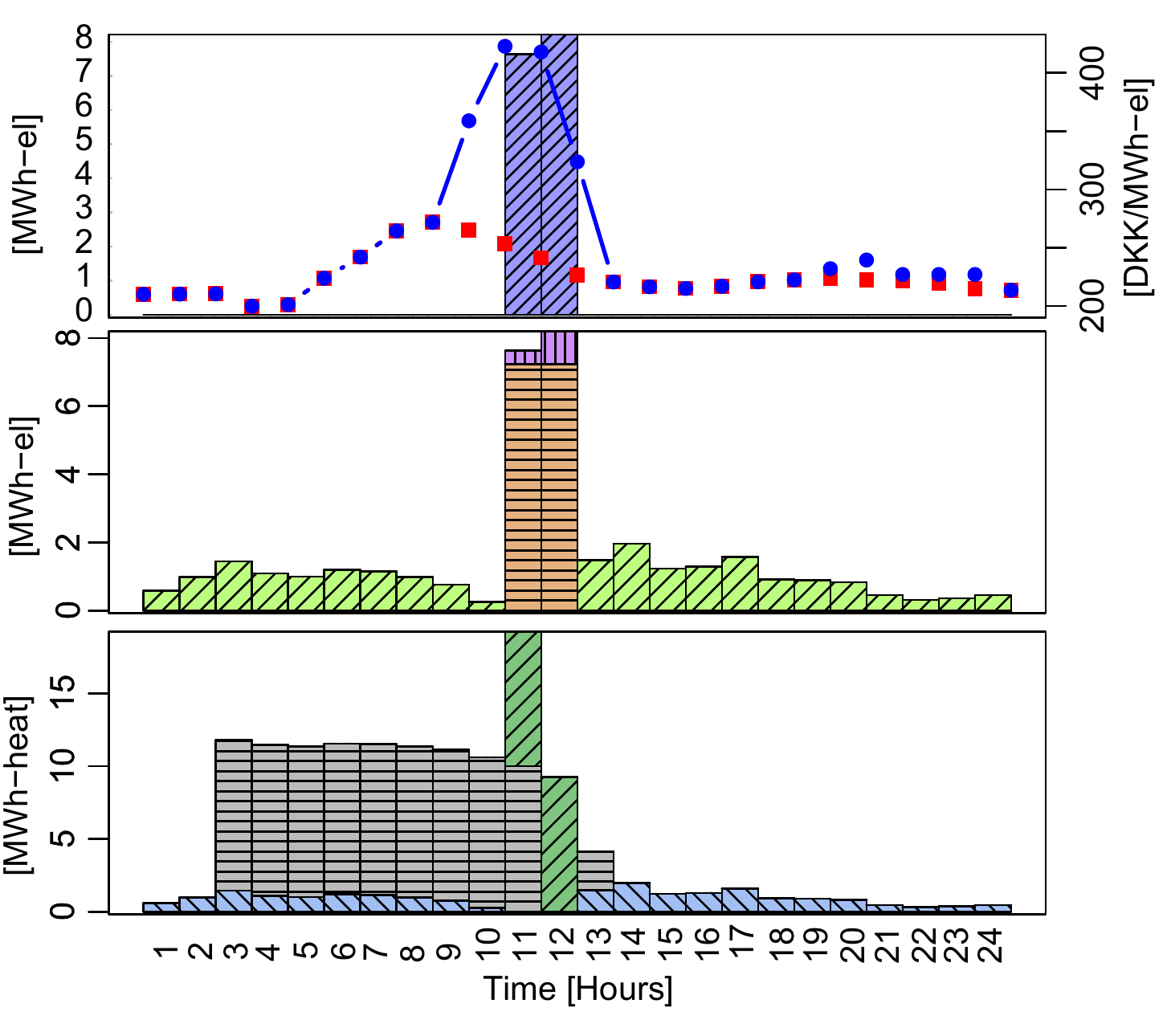}
	\caption{System operation: upward regulation amount including day-ahead and balancing prices (top), power production (middle), heat production (bottom)}
	\label{fig:Upcase2Operation}
	\end{subfigure}
	\caption{Upward regulation provided on 27th March 2017} \label{fig:UpCase2}
	
\end{figure}

	The second case for upward regulation is displayed in Figure \ref{fig:UpCase2} that follows the same structure as Figure \ref{fig:UpCase1}. Figure \ref{fig:Upcase2BiddingCurves} shows that two bids for upward regulation are won. As we can see in  Figure \ref{fig:Upcase2Operation} (middle), the upward balancing regulation is provided by the wind farm and the two CHP units in our system during these two hours. For this two hours the upward prices are significantly high and consequently, it is profitable to turn on the two CHP units. 

    Based on the system behaviour on those two representative days, we can summarize the two cases in which the DH operator can provide upward regulation. First, if we have an higher production of wind power than anticipated and offered in the day-ahead market. Second, if the upward regulation price is high enough that it is beneficial to start up the rather expensive CHP units.
  
\subsubsection{Downward regulation} \label{sec:downwardprices}

\begin{figure}
	\begin{subfigure}[t]{.45\textwidth}
	\centering
	\includegraphics[width=1.0\textwidth]{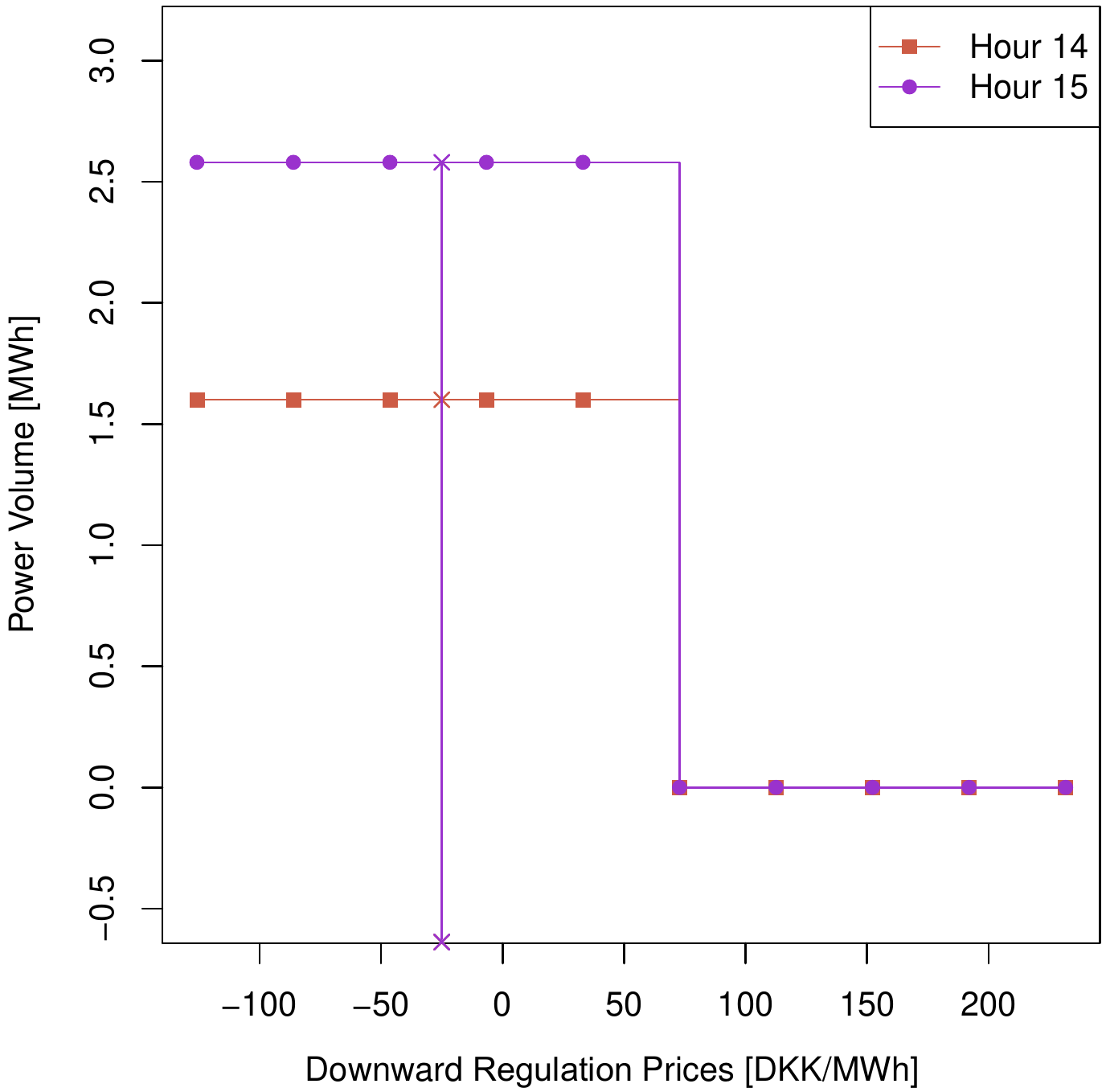}
	\caption{Bidding curves and won bids (marked with $\times$)}
	\label{fig:Downcase1BiddingCurves}
	\end{subfigure}
   \begin{subfigure}[t]{.52\textwidth}
	\centering
	\includegraphics[width=1.0\textwidth]{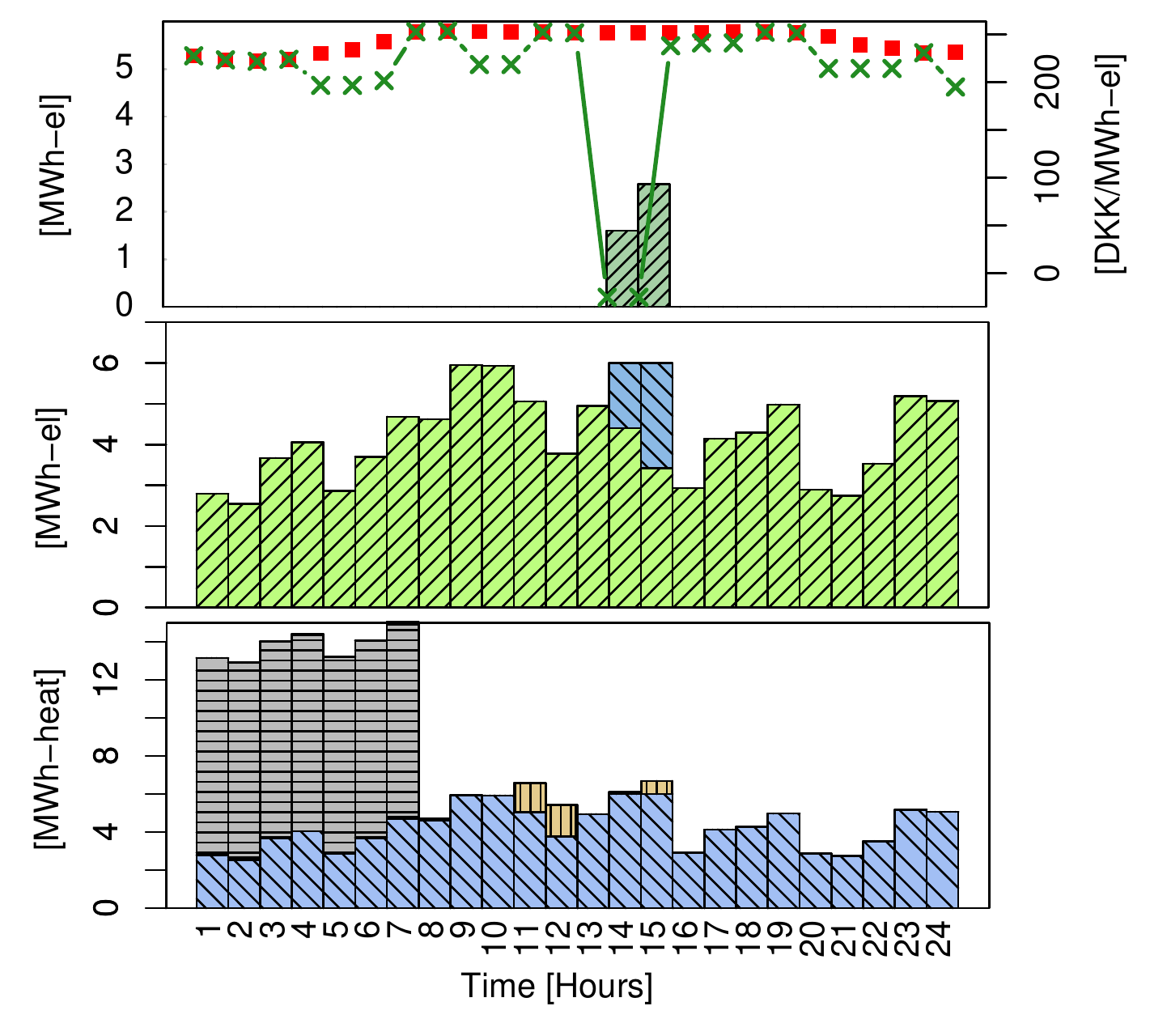}
	\caption{System operation: day-ahead committed amount, downward regulation amount as well as day-ahead and balancing prices (top), power production (middle), heat production (bottom)}
	\label{fig:Downcase1Operation}
	\end{subfigure}
	\caption{Downward regulation provided on 6th February 2017} \label{fig:DownCase1}
		\end{figure}

In the following we analyze how a DH operator can provide downward regulation. The first option is presented in Figure \ref{fig:DownCase1}, which shows downward regulation provided in hour 14 and 15. In this case, the model decides to buy electricity from the grid at the downward price and turn on the electric boiler (see Figure \ref{fig:Downcase1Operation} (middle). In general, electric boilers are good candidates to provide downward regulation because they can absorb large volumes of power in a very short time. Thus, producing heat using the electric boiler constitutes a very economical option when downward regulation is needed.   
		
\begin{figure}
		\begin{subfigure}[t]{.45\textwidth}
	\centering
	\includegraphics[width=1.0\textwidth]{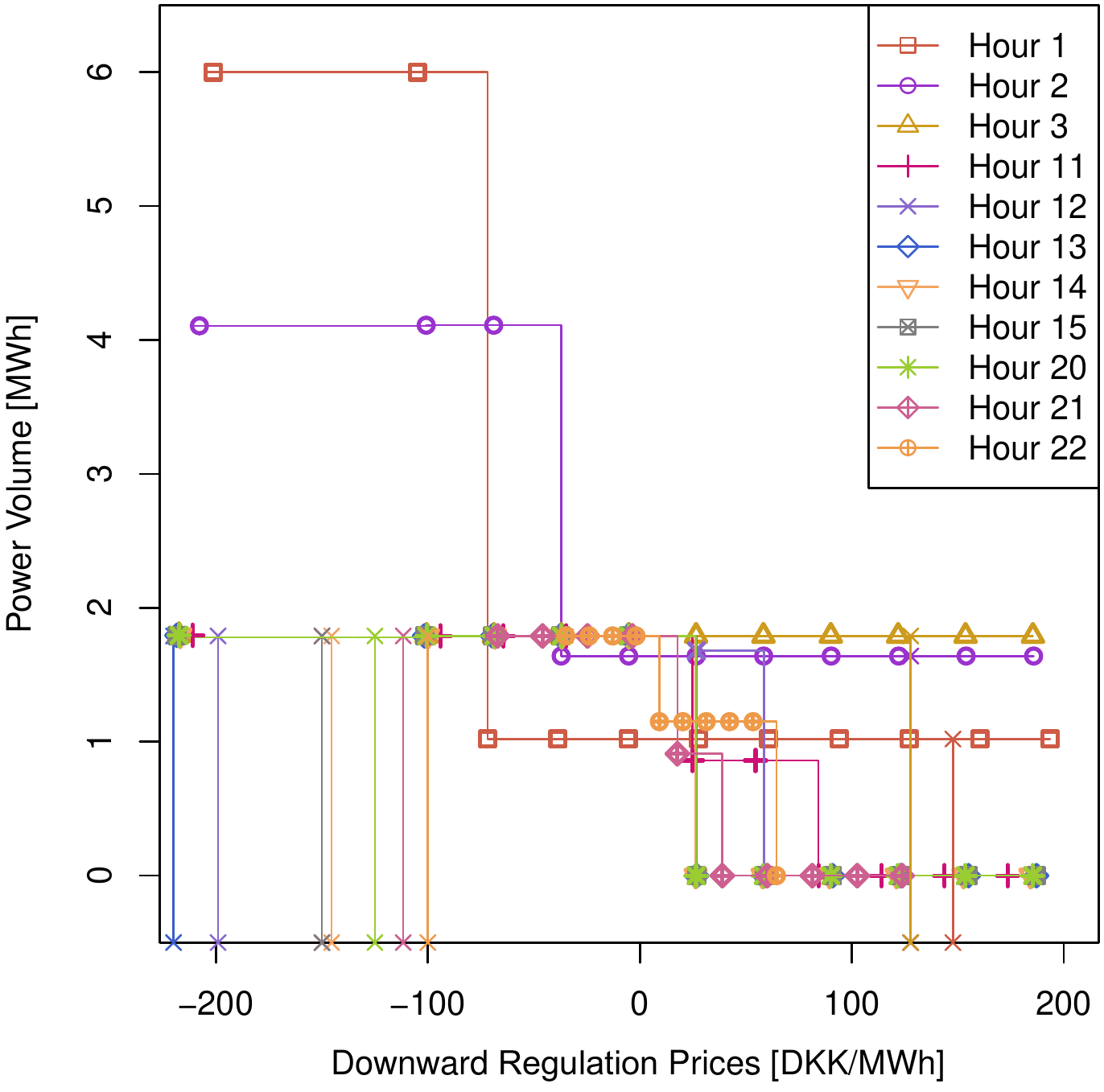}
	\caption{Bidding curves and won bids (marked with $\times$)}
	\label{fig:Downcase2BiddingCurves}
	\end{subfigure}
   \begin{subfigure}[t]{.52\textwidth}
	\centering
	\includegraphics[width=1.0\textwidth]{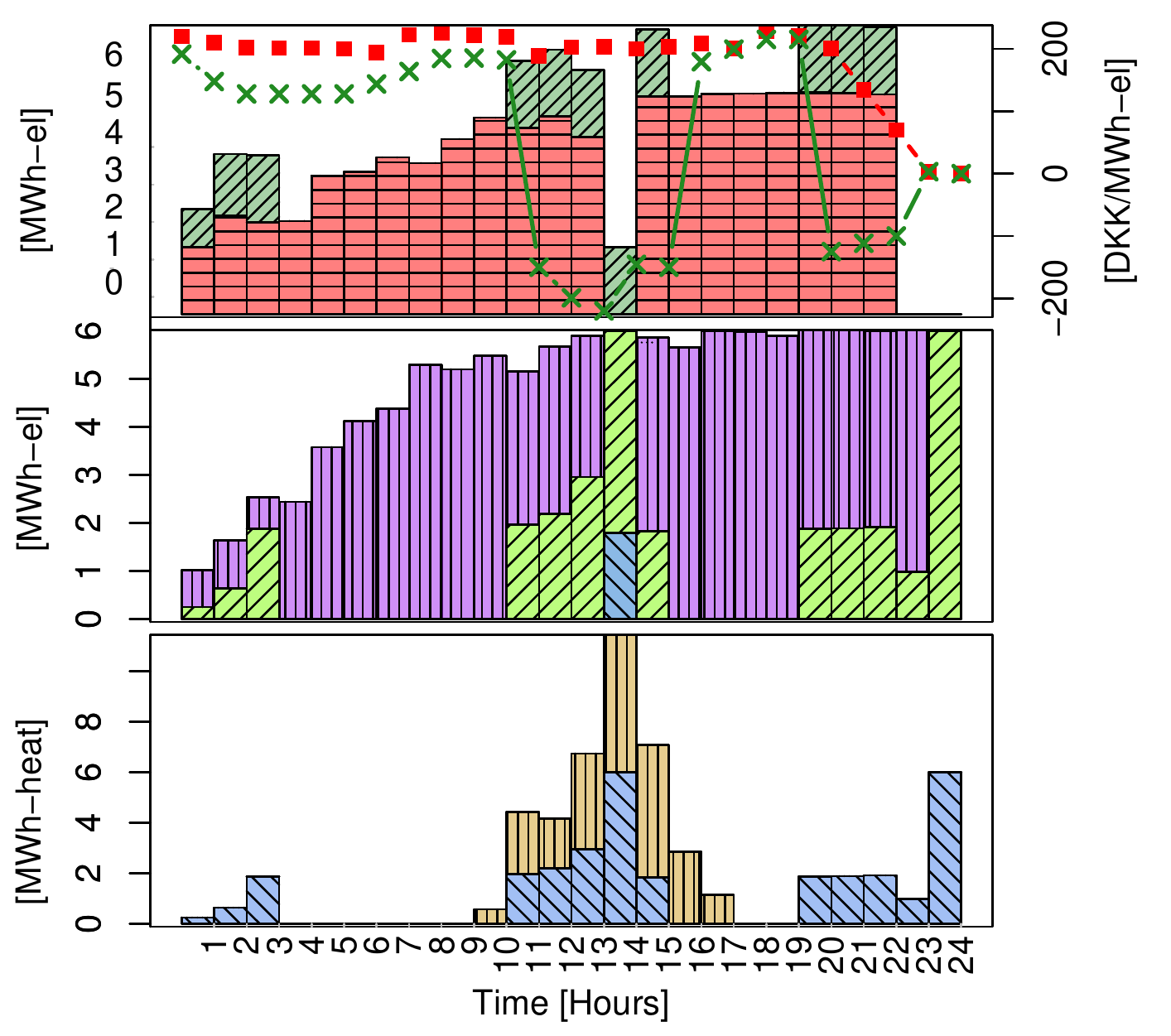}
	\caption{System operation: day-ahead committed amount, downward regulation amount as well as day-ahead and balancing prices (top), power production (middle), heat production (bottom)}
	\label{fig:Downcase2Operation}
	\end{subfigure}
	\caption{Downward regulation provided on 10th September 2017} \label{fig:DownCase2}
\end{figure}

	The second option in which our DH system can benefit from downward regulation is shown in Figure \ref{fig:DownCase2}. In this case the system takes advantage of the power sold previously on the day-ahead market to provide downward regulation. As it can be seen from Figure \ref{fig:Downcase2BiddingCurves} the system wins 11 bids for downward regulation on that specific day. Here the system stops providing the day-ahead power previously dispatched and buys this lack of production at the downward price. The profit of the system is the difference between the electricity sold at the day-ahead price and the electricity bought at the balancing price. This behaviour is shown in Figure \ref{fig:Downcase2Operation} (top), where the difference between the power sold in the day-ahead market and the one sold in the balancing market is the actual production of our wind generators sold to the market (Figure \ref{fig:Downcase2Operation} (middle)). This behaviour is the same for all time periods where downward regulation is provided with the exception of the hour 14 in which no day-ahead auction is won for that hour and therefore, the system decides to buy downward regulation and turn on the electric boiler (Figure \ref{fig:Downcase2Operation} (middle)).      
Based on the results in this section, we can summarize two ways of providing downward regulation for a DH operator. Either the electric boiler is used to provide downward regulation and produce at a low price or previously won power bids on the day-ahead market are corrected due to lower wind power production than excepted.

\section{Summary and outlook}

In this work, we present a planning method based on two-stage stochastic programming that allows DH providers, which operate a portfolio of units and have uncertain RES production, to create price dependent bids for both day-ahead and balancing markets and optimize the daily production. First, a stochastic program is solved to obtain and present the bids to the day-ahead market. Once the market is cleared, and the producer knows the power production plan, a second stochastic program is used on an hourly basis to generate bids for the balancing market considering the day-ahead power previously dispatched. After the bids for the balancing market are created and submitted, the market is cleared and the model optimizes the heat production for the new power commitment plan. In addition, we propose a new methodology to define balancing prices scenarios that account for the volatility of these latter based on their observed mean duration and values. 

We perform an extensive analysis of the production and trading behaviour of a real DH system in the two markets. The results confirm that uncertain electricity prices have a large impact on the system cost followed by uncertainty in the wind power production. In contrast, solar thermal production uncertainty has a minor influence due to the flexibility given by the heat tank storage. We also show the benefits of using a special tariff that utilizes the power production of wind farms with an electric boiler. This special tariff reduces the yearly total system cost enormously.
Regarding the inclusion of balancing market trading into the solution approach, we show that the participation in this market translates in larger profits resulting in lower operational costs. Finally, we investigate the behaviour of the system in case of upward and downward regulation in more detail. The results emphasize the important role of an electric boiler as flexible unit connected to the markets. 
To summarize, we propose a new planning method to reduce the impact of uncertainties on the production planning for DH systems. In order to achieve this, we hedge against uncertain electricity market prices and production using stochastic programming to create price dependent bids. The integration of RES production is facilitated by re-dispatching the imbalances in the balancing market. Furthermore, we show that considering the DH system as portfolio of units enables the necessary flexibility to react to seasonal changes and uncertainties.

We envision three different lines of future work. First, to use the presented approach to aggregate offers from a portfolio of different DH producers and calculate the optimally combined offer that can maximize the profit of all producers considering that we are now price-makers instead of price-takers. Therefore, a bi-level optimization program should be formulated. Second, to improve the bidding strategies to hedge even more against uncertain electricity prices. Finally, as our results indicate a significant margin of improvement by using the balancing market. Therefore, it becomes essential to develop more accurate forecasting techniques to predict balancing prices and their high volatility for one or two hours in advance. 
\label{sec:conclusion}


\subsubsection*{Acknowledgements}
This work is partly funded by Innovation Fund Denmark through the CITIES research center (no. 1035-00027B). We would like to thank Hvide Sande Fjernvarme A.m.b.A. especially Martin Halkjær Kristensen for the valuable input and data set.

\small

\bibliography{bibliography}

\end{document}